\documentclass[smallextended,envcountsect,]{svjour3}
\smartqed
\usepackage{graphicx}

\usepackage{latexsym, amsmath,amsxtra, amssymb, latexsym, amscd}
\usepackage{enumerate,cite,color}
\usepackage{hyperref}

\def\Limsup{\mathop{{\rm Lim}\,{\rm sup}}}
\def\Sol{\mbox{\rm Sol}\,}

\def\Limsup{\mathop{{\rm Lim}\,{\rm sup}}}

\def\epi{\mbox{\rm epi}\,}

\def\dom{\mbox{\rm dom}\,}

\def\R{\mathbb{R}}
\def\N{\mathbb{N}}
\def\Sol{\mbox{\rm Sol}}

\def\epi{\mathrm{epi}\,}

\def\dom{\mathrm{dom}\,}

\def\bd{\mathrm{bd}\,}

\begin{document}

\title{Directional  Subdifferentials at Infinity and Its Applications}

\titlerunning{Directional  Subdifferentials at Infinity and Its Applications}

\authorrunning{L.N. Kien, N.V. Tuyen, T.V. Nghi}


\author{Le Ngoc Kien  \and Nguyen Van Tuyen \and  Tran Van Nghi}

\institute{
           Le Ngoc Kien \at
           Faculty of Fundamental Science, Vietnam-Hungary Industrial University,
           \\
           Tung Thien, Hanoi, Vietnam;
           \\
           Department of Mathematics, Hanoi Pedagogical University 2,\\
           Xuan Hoa, Phu Tho, Vietnam \\
           lengockien@viu.edu.vn
           \and
           Nguyen Van Tuyen,  Corresponding author \at
           Department of Mathematics, Hanoi Pedagogical University 2,\\
           Xuan Hoa, Phu Tho, Vietnam \\
           nguyenvantuyen83@hpu2.edu.vn; tuyensp2@yahoo.com 
            \and 
          Tran Van Nghi, Co-corresponding author   \at
          Department of Mathematics, Hanoi Pedagogical University 2,\\
          Xuan Hoa, Phu Tho, Vietnam \\
          tranvannghi@hpu2.edu.vn              
}

\date{Received: date / Accepted: date}

\maketitle

\begin{abstract}
This paper investigates the behavior of sets and functions at infinity by introducing new concepts, namely directional normal cones at infinity for unbounded sets, along with  limiting and singular subdifferentials at infinity in the direction for extended real-valued functions. We develop several calculus rules for these concepts and then apply them to nonsmooth optimization problems. The applications include establishing directional optimality conditions at infinity, analyzing the coercivity, proving the compactness of the global solution set, and examining properties such as weak sharp minima and error bounds at infinity. To demonstrate the effectiveness of the proposed approach, illustrative examples are provided and compared with existing results. 
\end{abstract}


\section{Introduction}\label{Introduction}

The Mordukhovich (limiting) normal cone and its associated limiting subdifferential play a key role in modern variational analysis and optimization, as they provide useful tools for addressing nonsmooth and nonconvex optimization problems (see \cite{Mordukhovich2006,Mordukhovich2018,Mordukhovich2024}). 
\medskip

The study of the behavior of sets and functions at infinity is a significant topic in optimization theory. Sufficient conditions for the existence of solutions in nonsmooth optimization via asymptotic cones and generalized asymptotic functions were given in \cite{AT03,Aus1,Aus2,FFB-Vera,F-L-Vera,HLL,HLM,KMPT-21,Lara-Tuyen-Nghi,Nghi-Kien-Tuyen}. Recently, the concepts of normal cones at infinity for unbounded sets, together with limiting and singular subdifferentials at infinity for extended real-valued functions, were introduced in \cite{Kim-Tung-Son-23}. By employing various calculus rules for these notions, the authors characterized Lipschitz continuity at infinity for lower semicontinuous functions and applied the results to optimization problems, including optimality conditions, weak sharp minima, and stability properties. In \cite{Tung-Son-24}, Clarke’s tangent cones at infinity for unbounded sets, subgradients at infinity for extended real-valued functions, and necessary optimality conditions at infinity for optimization problems were studied. Sufficient conditions for the existence of error bounds at infinity in lower semicontinuous inequality systems, as well as necessary optimality conditions for constrained optimization problems, were established in \cite{Tuyen-24}. By using the tool of subdifferentials at infinity, Tuyen, Bae, and Kim  \cite{Tuyen-25} proposed optimality conditions at infinity for nonsmooth minimax programming problems. 

Very recently, Anh and Hung   \cite{HA-Hung-25} investigated properties of normal cones with respect to a set and developed calculus rules for subdifferentials relative to a set at a reference point. The obtained results extend and improve the corresponding ones in \cite{Thinh-Chuong-18,Thinh-Chuong-19,Thinh-Qin-Yao}. Furthermore, the authors also introduced the notions of normal cones and subdifferentials with respect to a set at infinity. These tools were then employed to derive necessary optimality conditions at infinity, establish the compactness of the solution set, and verify the coercivity  in optimization problems with unbounded feasible sets. It is well known that, in the theory of subdifferentials (or generalized derivatives), both the size of the subdifferential and the availability of computation rules play a crucial role. In \cite{HA-Hung-25}, the authors provided examples showing that the subdifferentials at infinity introduced in \cite{Kim-Tung-Son-23,Tung-Son-24} can be quite large.  However, for the subdifferential with respect to a set at infinity, as introduced in \cite{HA-Hung-25}, the computation is quite involved, and the fundamental formulas for calculation have not yet been fully established. Therefore, the development of new subdifferential concepts at infinity remains an important direction of research.

In this paper, motivated by \cite{Ginchev-Mor-11,Ginchev-Mor-12,Kim-Tung-Son-23}, we study    normal cone at infinity and subdifferentials at infinity in the direction and their applications.  Our main contributions to variational analysis and nonsmooth optimization theory are summarized as follows:
\medskip

$\bullet$\, For variational analysis, we introduce new concepts: the directional normal cone at infinity for unbounded sets, and the directional limiting and the directional singular  subdifferentials at infinity  for extended real-valued functions. Several fundamental calculus rules for these notions are established. Furthermore, we prove that the directional Lipschitz property at infinity of a real-valued function is equivalent to that  its directional singular subdifferential at infinity is equal to $\{0\}$.
\medskip

$\bullet$\, For nonsmooth optimization theory, by employing related properties and calculus rules, we present several applications to nonsmooth optimization problems, including directional optimality conditions at infinity, coercivity,  compactness of the global solution set, the weak sharp minima property at infinity, and the error bound property at infinity.

This paper is organized as follows. Section \ref{Preliminaries} reviews some necessary definitions and preliminary results from variational analysis. In Section \ref{Section-3}, we introduce and study the directional normal cone at infinity, along with the directional  limiting and singular subdifferentials at infinity. Section \ref{Section-4} is devoted to applications in nonsmooth optimization problems. Finally, Section \ref{Conclusion} presents a discussion of problems for future research.

\section{Preliminaries}\label{Preliminaries}
In this section, we recall several notions related to generalized differentiation from \cite{AT03,Kim-Tung-Son-23,Mordukhovich2006,Mordukhovich2018,Penot-2013,Rockafellar1998}. 

\subsection{Notation and Definition} 
Throughout the paper, denote $\N:=\{1, 2, \ldots\}$ and let $\R^n$ be the Euclidean space with the usual scalar product $\langle \cdot, \cdot\rangle$ 
the corresponding  Euclidean norm $\|\cdot\|$ where $n\in \N$.  The closed unit ball   and the nonnegative orthant  in $\R^n$ are denoted, respectively, by $\mathbb{B}$ and $\R^n_+$. The closed ball centered at the origin with radius $R>0$ is denoted by $\mathbb{B}_R$. Let $D$  be a subset of $\R^n$. We say that $D$ is {\em locally closed}  if for any $x\in D$ there is a neighborhood $U$ of $x$ such that $D\cap U$ is closed. The interior, the boundary, and the convex hull of  $D$ are  denoted, respectively, by $\mathrm{int}\,D$, $\mathrm{bd}\,D$, and  $\mathrm{co}\, D$. For a given point $x\in\R^n$, we denote the  Euclidean projector of $x$ onto $D$ and the distance from $x$ to $D$  by $\Pi_D(x)$ and $\mathrm{dist}(x; D)$ (or, $d_D(x)$), respectively.     As usual,   $\mathbb{S}:=\{x\in \R^n \mid \|x\|=1\}$ is the unit sphere in $\R^n$,  $\overline{\mathbb{R}}:=\R\cup\{\infty\}$ is the extended real line,  $[\alpha]_+:=\max\,\{\alpha, 0\}$ for any $\alpha\in\R$. The notation $x\to \infty$ means that $\|x\|\to \infty$. The {\em asymptotic cone} 
of $D$, denoted by $D^\infty$, is  defined by
$$D^{\infty}:=\left \{  u\in \mathbb{R}^{n}\,|\,\exists~t_{k}\rightarrow +
\infty,  x_{k}\in D,~\frac{x_{k}}{t_{k}}\rightarrow u\right \}.$$
When $D=\emptyset$, we put $\emptyset^{\infty}:=\emptyset$.    It follows from \cite[Proposition 2.1.2]{AT03} that $D$ is {\em  bounded} if and only if $D^\infty=\{0\}$. 

\medskip
Given an extended real-valued function $f: \mathbb{R}^n\rightarrow \overline{\mathbb{R}}$. The \textit{effective domain} and the \textit{epigraph} of $f$ are denoted, respectively, by
$$\mbox{dom}\, f:=\{x \in \mathbb{R}^n \mid f(x) < +\infty\}$$ 
and  
$$ {\rm{epi}}\, f:=\{ (x, \alpha) \in \mathbb{R}^n \times \mathbb{R} \mid \alpha \ge f (x)\}.$$  
We say that $f$ is {\em proper} (resp., {\em proper at infinity}) if  $\mathrm{dom} f$ is  nonempty (resp., $\mathrm{dom} f$ is unbounded).  The function $f$ is called {\em  lower semicontinuous} (l.s.c.) if its epigraph is closed.

Let  $F : \R^n \rightrightarrows \R^m$ be a set-valued mapping. The \textit{domain} and the \textit{graph}  of $F$ are given, respectively, by
$${\rm dom}\,F:=\{x\in \R^n \mid F(x)\not= \emptyset\} $$
and
$$ {\rm gph}\,F:=\{(x,y)\in \R^n \times \R^m \mid y \in F(x)\}.$$
The set-valued mapping $F$ is called   \textit{proper} if $\dom F \not= \emptyset.$

\medskip
For a set-valued map $F : \R^n \rightrightarrows \R^m$, the \textit{Painlev\'e-Kuratowski outer/upper limit} of $F$ at $\bar x$ is defined by
\begin{align*} 
	\Limsup\limits_{x\rightarrow \bar x} F(x):=\Big\{ y\in \mathbb{R}^m \mid \exists x_k \rightarrow \bar x, y_k \rightarrow y \ \mbox{with}\ y_k\in F(x_k) \ \ \forall k\in\N\Big\}
\end{align*}
and as $x\to \infty$ we define
\begin{align*} 
	\Limsup\limits_{x\rightarrow \infty} F(x):=\Big\{ y\in \mathbb{R}^m \mid \exists x_k \rightarrow \infty, y_k \rightarrow y \ \mbox{with}\ y_k\in F(x_k)\ \ \forall k\in\N\Big\}.
\end{align*}

\subsection{Normal Cones and Subdifferentials} 
We summarize in this subsection several notions from variational analysis, in particular the concepts of normal cones and subdifferentials, following  \cite{Mordukhovich2006,Mordukhovich2018}.

\begin{definition}{\rm 
		Let $\Omega$ be a nonempty subset of $\mathbb{R}^n$ and let $\bar x \in \Omega$. 
		\begin{enumerate}[(i)]
			\item The \textit{regular/Fr\'echet normal cone} to $\Omega$ at $\bar x$ is defined by
			\begin{align*}
				\widehat N(\bar x; \Omega)=\left\{ v\in \mathbb{R}^n\mid \limsup\limits_{x \xrightarrow{\Omega}\bar x} \dfrac{\langle v, x-\bar x \rangle}{\|x-\bar x\|} \leq 0 \right\},
			\end{align*}
			where $x \xrightarrow{\Omega} \bar x$ means that $x \rightarrow \bar x$ and $ x\in \Omega$.
			\item The \textit{limiting/Mordukhovich normal cone} to $\Omega$ at $\bar x$ is given by
			\begin{align*}
				N(\bar x; \Omega)=\Limsup\limits_{ x \xrightarrow{\Omega} \bar x} \widehat{N}(x; \Omega).
			\end{align*}
			When $\bar x \not\in \Omega$, we put $\widehat N(\bar x;\Omega)=N(\bar x;\Omega) :=\emptyset$.
	\end{enumerate} }
\end{definition}	
By the above definitions, it is clear that
\begin{align*} 
	\widehat N(x;\Omega) \subset N(x;\Omega) \ \    \forall  x \in \Omega.
\end{align*}

\begin{definition}\label{def21} {\rm Given a function   $f \colon \mathbb{R}^n \to \overline{\mathbb{R}}$ and a point  $\bar{x} \in {\rm dom}f$.  
		\begin{enumerate}[{\rm (i)}]
			\item The {\em regular/Fr\'echet subdifferential} of $f$ at $\bar{x}$ is defined by 
			$$
			\widehat{\partial}f(\bar{x}):=\{ v \in \mathbb{R}^n \mid (v,-1)\in \widehat{N}((\bar{x},f(\bar{x}));\mathrm{epi} f)  \}.  
			$$
			\item The {\em limiting/Mordukhovich subdifferential} and the {\em limiting/Mordukhovich singular subdifferential}  of $f$ at $\bar{x}$ are defined, respectively, by  
			$$
			\partial f(\bar{x}):=\{ v \in \mathbb{R}^n \mid (v,-1)\in {N}((\bar{x},f(\bar{x}));\mathrm{epi} f)  \}, 
			$$
			and   
			$$\partial^{\infty} f(\bar{x}):=\{ v \in \mathbb{R}^n \mid \exists v_k \in 	\widehat{\partial} f(\bar{x}), \lambda_k \downarrow 0, \lambda_k v_k \to v \}.	$$
	\end{enumerate}}
\end{definition}
It is well-known that
$$
\partial f(\bar{x})=\Limsup_{x \xrightarrow{f} \bar{x}}\widehat{\partial}f(x) \supseteq \widehat{\partial}f(x),
$$
where $x \xrightarrow{f} \bar{x}$ means that $x \to \bar{x}$ and $f(x) \to f(\bar{x})$. For a convex function  $f$, the subdifferentials $\widehat{\partial} f(\bar{x})$ and $\partial f(\bar{x})$ coincide with the subdifferential in the sense of convex analysis. We obtain that
\begin{align*}
	\partial^{\infty} f(\bar{x})\subseteq \{ v \in \mathbb{R}^n \mid (v,0)\in N((\bar{x},f(\bar{x}));\mathrm{epi} f)  \},
\end{align*}
and the inclusion holds with equality whenever $f$ is locally l.s.c. at $\bar{x}$ (see \cite[Theorem 8.9]{Rockafellar1998}). 

Given a set $\Omega \subset \mathbb{R}^n$. The {\em indicator function} $\delta_{\Omega} \colon \mathbb{R}^n \to \overline{\mathbb{R}}$ of $\Omega$ is defined by
$$\delta_{\Omega}(x) :=
\begin{cases}
	0 & \textrm{ if } x \in \Omega, \\
	+\infty & \textrm{ otherwise.}
\end{cases}$$
For any $x \in \Omega$, it holds that $\partial \delta_{\Omega} (x) =\partial^{\infty} \delta_{\Omega} (x)= N(x; \Omega)$ (see \cite{Mordukhovich2006,Mordukhovich2018,Rockafellar1998}).

A nonsmooth versions of Fermat's rule is stated as follows.
\begin{lemma}[{see \cite[Proposition 1.114]{Mordukhovich2006}}] \label{lema22}
	If a proper function $f \colon \mathbb{R}^n \to \overline{\mathbb{R}}$ has a local minimum at $\bar{x},$ then $0 \in \widehat{\partial}f(\bar{x})\subset \partial f(\bar{x}).$	
\end{lemma}

The following  Ekeland variational principle  is an useful tool in establishing our main results.
\begin{lemma}[{see \cite[Theorem 1]{Ekeland-74}}]\label{lema26} 
	Let $f \colon \mathbb{R}^n \to \overline{\mathbb{R}}$ be a proper l.s.c. function and bounded from below. Let $\epsilon >0$ and $u \in \mathbb{R}^n$ be satisfied
	\begin{eqnarray*}
		f (u) &\le& \inf_{x \in \mathbb{R}^n}f (x) + \epsilon. 
	\end{eqnarray*}
	Then, for any $\lambda >0$ there exists  $v \in \mathbb{R}^n$ such that
	\begin{enumerate}[{\rm (i)}]
		\item $f (v) \le f(u),$
		\item $\|v-u\| \le \lambda,$ and
		\item $f (v) \le  f(x)+\dfrac{\epsilon}{\lambda}\|x-v\|$ for all $x \in \mathbb{R}^n.$
	\end{enumerate}
\end{lemma}

\subsection{Normal Cones and Subdifferentials at Infinity}
This subsection provides a brief review of the notions of normal cones and subdifferentials at infinity, as presented in \cite{Kim-Tung-Son-23}. 

\begin{definition}\label{def31} {\rm 
		Let $\Omega$ be an unbounded subset in $\R^n$. The {\em norm cone to the set $\Omega$ at infinity} is defined by
		\begin{eqnarray*}
			N(\infty; \Omega) &:=& \Limsup_{ x \xrightarrow{\Omega} \infty} \widehat{N}(x; \Omega),
		\end{eqnarray*} 
		where $x \xrightarrow{\Omega} \infty$ means that $\|x\|\to \infty$ and $x\in \Omega$.
}\end{definition}
By \cite[Propositions 3.5 and 3.6]{Kim-Tung-Son-23}, we have 
\begin{eqnarray*}
	N(\infty; \Omega) &:=& \Limsup_{ x \xrightarrow{\Omega} \infty} {N}(x; \Omega),
\end{eqnarray*}
and  $N(\infty; \Omega)$ is nontrivial  if and only if $\mathrm{bd}\,\Omega$ is unbounded. Furthermore, by definition, it is easy to see that
\begin{equation*}
	N(\infty; \Omega)=\bigcup_{i=1}^n N(\infty_{\{i\}}; \Omega),
\end{equation*}
where 
\begin{equation*}
	N(\infty_{\{i\}}; \Omega) :=\Limsup_{ x\in \Omega, |x_i|\to\infty} \widehat{N}(x; \Omega), \ \ i=1, \ldots, n.
\end{equation*}

\begin{definition}\label{def41} {\rm 
		Let $f\colon\R^n\to \overline{\R}$ be an l.s.c. and proper at infinity function. The {\em limiting/Mordukhovich and the singular subdifferentials} of $f$ at infinity are defined, respectively, by 
		\begin{eqnarray*}
			\partial f(\infty) &:=& \{u \in \mathbb{R}^n \ | \ (u, -1) \in \mathcal{N} \},\\
			\partial^{\infty} f(\infty) &:=& \{u \in \mathbb{R}^n \ | \ (u, 0) \in \mathcal{N}\},
		\end{eqnarray*}
		where $\mathcal{N} := \displaystyle \Limsup_{x \to \infty} {N}((x,f(x)); \textrm{epi} f).$
}\end{definition}

The next result provides limiting characterizations of both limiting and singular subdifferentials at infinity.
\begin{proposition}[{see \cite[Proposition 4.4]{Kim-Tung-Son-23}}]\label{pro42} 
	The following relationships hold
	\begin{equation*}
		\aligned
		\partial f(\infty) &= \Limsup_{x \to \infty} \partial f(x)= \Limsup_{x \to \infty} \widehat{\partial} f(x), 
		\\
		\partial^{\infty} f(\infty) &= \Limsup_{x \to \infty, r \downarrow 0} r \partial f(x)\supseteq \Limsup_{x \to \infty} \partial^{\infty} f(x).  
		\endaligned
	\end{equation*}
\end{proposition}

For an unbounded closed set $\Omega\subset \R^n$,  it follows from Proposition \ref{pro42} and \cite[Proposition 1.19]{Mordukhovich2018} that
\begin{equation*} 
	\partial\delta_\Omega(\infty)=\partial^\infty\delta_\Omega(\infty)=N(\infty; \Omega).
\end{equation*}

Let us recall the definition of the Lipschitz property at infinity for l.s.c. functions, as introduced in \cite{Kim-Tung-Son-23,Tung-Son-24}. 
\begin{definition}\label{def51}{\rm
		Let  $f \colon \mathbb{R}^n \to \mathbb{R}$ be an l.s.c. function.  We say that $f$ is {\em Lipschitz at infinity} if there exist constants $L > 0$ and $R > 0$ such that
		\begin{eqnarray*}
			|f(x) - f(x')| &\le& L \|x - x'\| \ \  \textrm{ for all } \ \ x, x' \in \mathbb{R}^n \setminus \mathbb{B}_R.
		\end{eqnarray*}
}\end{definition}
The following theorem establishes a necessary and sufficient condition for the Lipschitz property at infinity of l.s.c. functions. 
\begin{proposition}[{see \cite[Proposition 5.2]{Kim-Tung-Son-23}}]\label{pro52}   
	Let $f \colon \mathbb{R}^n \to \mathbb{R}$ be an l.s.c. function. Then $f$ is Lipschitz at infinity if and only if $\partial^{\infty}f(\infty)=\{0\}.$ In this case, $\partial f(\infty)$ is nonempty compact. 
\end{proposition}
\section{Directional Normal Cone and Subdifferential at Infinity}\label{Section-3} 
In this section, we define and study directional normal cones at infinity for unbounded sets, along with directional limiting and singular subdifferentials at infinity for extended real-valued functions.
\subsection{Directional Normal Cone at Infinity}
Let $\Omega$ be a locally closed and unbounded  subset in $\R^n$ and $u\in\mathbb{S}$. The directional normal cone at infinity is given by the following definition.
\begin{definition}\rm   The {\em  normal cone to $\Omega$ in direction $u$ at infinity}, denoted by $N_\Omega(\infty; u)$, is defined by
	\begin{equation*}
		N_\Omega(\infty; u):=\Limsup_{x \xrightarrow{\Omega, u} \infty} \widehat{N}_\Omega(x),
	\end{equation*}
	i.e.,
	\begin{equation*}
		N_\Omega(\infty; u):=\left\{\xi\in\R^n\,|\, \exists x_k\xrightarrow{\Omega}\infty, \tfrac{x_k}{\|x_k\|}\to u, \xi_k\in \widehat{N}_\Omega(x_k), \xi_k\to \xi\right\}.
	\end{equation*}
\end{definition}
By definition, if $u\notin\Omega^\infty$, then $N_\Omega(\infty; u)=\emptyset$. Hence, we only consider the case that $u\in \Omega^\infty\cap\mathbb{S}$.

\medskip
To highlight the difference between the directional normal cone at infinity and the normal cone at infinity (see \cite{Kim-Tung-Son-23}), we present the example below.
\begin{example}\rm  Let $\Omega\subset\R^2$ be the set defined by
	\begin{equation*}
		\Omega:=\left\{x=(x_1, x_2)\in\R^2\\,|\, \tfrac{x_1}{2}\leq x_2\leq 2x_1, x_1\geq 0\right\}.
	\end{equation*}
	Let $u=(\tfrac{2}{\sqrt{5}},\tfrac{1}{\sqrt{5}})$ and $v=(\tfrac{1}{\sqrt{5}},\tfrac{2}{\sqrt{5}})$. Then, we have $u, v\in\Omega^\infty\cap\mathbb{S}$, and 
	\begin{align*}
		N_\Omega(\infty; u)&=\{(x_1, x_2)\,|\, x_2=-2x_1, x_1\geq 0\},
		\\
		N_\Omega(\infty; v)&=\{(x_1, x_2)\,|\, x_2=-\tfrac{1}{2}x_1, x_1\leq 0\},
	\end{align*}
	while
	\begin{equation*}
		N_\Omega(\infty)=N_\Omega (\infty_{\{1\}})=N_\Omega(\infty_{\{2\}})=N_\Omega(\infty; u)\cup N_\Omega(\infty; v).
	\end{equation*}
\end{example}

The relation between the directional normal cone at infinity and the normal cone at infinity is stated below.
\begin{proposition}\label{Proposition-3.1} 
	The following relation holds
	\begin{equation*}
		N_\Omega(\infty)=\bigcup_{u\in\Omega^\infty\cap\mathbb{S}}N_\Omega(\infty; u).
	\end{equation*}
\end{proposition}
\noindent{\it Proof\,}  The inclusion ``$\supset$'' is directly from the definition. Now, we let any $\xi\in N_\Omega(\infty)$. Then by definition, there exist sequences $x_k\xrightarrow{\Omega}\infty$ and $\xi_k\in\widehat{N}_\Omega(x_k)$ with $\xi_k\to\xi$ as $k\to\infty$. Without loss of generality, we may assume that $x_k\neq 0$ for all $k\in\N$. Consequently, the sequence $\tfrac{x_k}{\|x_k\|}$ converges to some $u\in \mathbb{S}$, taking a subsequence if necessary. Hence, $u\in\Omega^\infty$ and $\xi\in N_\Omega(\infty; u)$, as required. \qed

\medskip
The proposition below offers an essential method to determine $N_\Omega(\infty; u)$.
\begin{proposition} \label{pro1} For any $u\in\Omega^\infty\cap\mathbb{S}$, we have
	\begin{equation*}
		N_\Omega(\infty; u):=\Limsup_{x \xrightarrow{\Omega, u} \infty}{N}_\Omega(x).
	\end{equation*}
\end{proposition}
\noindent{\it Proof\,} Since $\widehat{N}_\Omega(x)\subset N_\Omega(x)$ for all $x\in\R^n$, one has
	\begin{equation*}
		N_\Omega(\infty; u)=\Limsup_{x \xrightarrow{\Omega, u} \infty} \widehat{N}_\Omega(x)\subset \Limsup_{x \xrightarrow{\Omega, u} \infty} {N}_\Omega(x).
	\end{equation*}
	We now take any $\xi\in \Limsup_{x \xrightarrow{\Omega, u} \infty} {N}_\Omega(x)$. By definition, there exist sequences $x_k\xrightarrow{\Omega}\infty$, $\xi_k\in N_\Omega(x_k)$ with $\tfrac{x_k}{\|x_k\|}\to u$  and $\xi_k\to \xi$. Hence, for each $k\in\N$, there are $z_k\in\Omega$ and $\eta_k\in\widehat{N}_\Omega(z_k)$ such that $\|z_k-x_k\|\leq \tfrac{1}{k}$ and $\|\eta_k-\xi_k\|\leq \tfrac{1}{k}$. Clearly, $z_k\to\infty$ and $\eta_k\to \xi$. We claim that $\tfrac{z_k}{\|z_k\|}\to u$ as $k\to\infty$. Indeed, since $z_k\to\infty$ and $x_k\to\infty$, we have
	\begin{align*}
		\left\|\frac{z_k}{\|z_k\|}-\frac{x_k}{\|x_k\|}\right\|&=\frac{\big\|z_k\|x_k\|-x_k\|z_k\|\big\|}{\|z_k\|\|x_k\|}\\
        &\leq  \frac{\big\|z_k\|x_k\|-x_k\|z_k\|\big\|}{\|z_k\|}
		\\
		&=\left\|z_k\frac{\|x_k\|}{\|z_k\|}-x_k\right\|
		\\
		&\leq \left\|z_k\frac{\|x_k\|}{\|z_k\|}-z_k\right\| +\|z_k-x_k\|
		\\
		&=|\|x_k\|-\|z_k\||+\|z_k-x_k\|\leq 2 \|z_k-x_k\|.
	\end{align*}
	This implies that
	\begin{equation*}
		\lim_{k\to\infty}\left\|\frac{z_k}{\|z_k\|}-\frac{x_k}{\|x_k\|}\right\|=0
	\end{equation*}
	and so $\tfrac{z_k}{\|z_k\|}\to u$ as $k\to\infty$. Hence, $\xi\in N_\Omega(\infty; u)$. The proof is complete. \qed
 
\medskip
The next proposition characterizes the nontriviality of $N_\Omega(\infty; u)$.
\begin{proposition} \label{pro36} Let $u\in \Omega^\infty\cap  \mathbb{S}$. Then,  $N_{\Omega}(\infty; u)\neq\{0\}$ if and only if  $u\in (\bd\Omega)^\infty$.
\end{proposition}
\noindent{\it Proof\,} Suppose that there exists  $\xi \in N_{\Omega}(\infty; u)\setminus\{0\}$. Then, by Proposition \ref{pro1}, there exist  $x_k\xrightarrow{\Omega, u}\infty$ and $\xi_k\in  N_{\Omega}(x_k)$ such that $\xi_k\to \xi$. This implies that  $\xi_k\neq0$ for all  $k$ large enough. By \cite[Proposition 1.2]{Mordukhovich2018}, $x_k\in \bd \Omega$ for all  $k$ large enough. This and the fact that $x_k\xrightarrow{\Omega, u}\infty$ imply that $u\in (\bd\Omega)^\infty$, as required.
	
	We now assume that $u\in (\bd\Omega)^\infty\cap \mathbb{S}$. By definition, there exist sequences $x_k\in\bd\Omega$ and $t_k\to \infty$ such that $\tfrac{x_k}{t_k}\to u$ with $\|u\|=1$. Hence, $\tfrac{\|x_k\|}{t_k}\to 1$ as $k\to\infty$. This implies that
	$\|x_k\|=t_k.\tfrac{\|x_k\|}{t_k}\to \infty$ and 
	$$\frac{x_k}{\|x_k\|}=\frac{x_k}{t_k}.\frac{t_k}{\|x_k\|}\to u.$$
	Hence, $x_k\xrightarrow{\bd\Omega, u}\infty$. Thus, by \cite[Proposition 1.2]{Mordukhovich2018}, for each $k$  there exists $\xi_k\in N_{\Omega}(x_k)$ with $\|\xi_k\|=1$. Without loss of generality, we may assume  that $\xi_{k}\to \xi$ for some $\xi\in \R^n$ with $\|\xi\|=1$. Since $x_k\in\bd\Omega$ and $\Omega$ is locally closed, it is clear that $x_k\in \Omega$ for all $k\in\N$. Hence $x_k\xrightarrow{\Omega, u}\infty$.  By Proposition \ref{pro1},  $\xi\in N_{\Omega}(\infty; u).$ The proof is complete.\qed

\medskip
The concept of a directional neighborhood at infinity is introduced as follows.
\begin{definition}\rm  Let $u\in\R^n$. A {\em neighborhood of the infinity in direction $u$} is defined by
	\begin{equation*}
		V_{R,\delta}(\infty; u):=\Big\{z\in\R^n\setminus\mathbb{B}_R\,|\, \big\|z\|u\| -u\|z\|\big\|\leq \delta\|z\|\|u\|\Big\}
	\end{equation*}
	for some $R>0$ and $\delta>0$.
\end{definition}

By definition, one has
\begin{equation*}
	V_{R,\delta}(\infty; u)=
	\begin{cases}
		\R^n\setminus \mathbb{B}_R,\ \ &\text{if}\ \ u= 0,
		\\
		\left\{z\in \R^n\setminus \mathbb{B}_R\,|\, \left\|\frac{z}{\|z\|}-\frac{u}{\|u\|}\right\|\leq \delta\right\},\ \ &\text{if}\ \ u\neq 0.
	\end{cases}
\end{equation*}

The following result provides the computation rule for the directional normal cone at infinity of the intersection of two sets.
\begin{proposition} \label{pro37} Let $\Omega_1$, $\Omega_2$ be unbounded subsets in $\R^n$. Then, for each $u\in (\Omega_1\cap\Omega_2)^\infty \cap \mathbb{S}$, if  the following condition is satisfied 
	\begin{equation}\label{37a}
		N_{\Omega_1}(\infty; u)\cap (-N_{\Omega_2}(\infty; u))=\{0\},
	\end{equation}
	then
	\begin{equation}\label{37b}
		N_{\Omega_1\cap \Omega_2}(\infty; u)\subset  N_{\Omega_1}(\infty; u)+N_{\Omega_2}(\infty; u).
	\end{equation}
\end{proposition}
\noindent{\it Proof\,} 	Firstly, we show  that there exists a constant  $R>0$  satisfying
	\begin{equation}\label{37c}
		N_{\Omega_1}(x)\cap (-N_{\Omega_2}(x))=\{0\} \; \forall x\in \Omega_1\cap \Omega_2\cap V_{R,\frac{1}{R}}(\infty; u).
	\end{equation}
	Indeed, on the contrary, suppose that for each $k$ large enough, there  exist sequence $x_k\in \Omega_1\cap \Omega_2\cap V_{k,\frac{1}{k}}(\infty; u)$ and $\xi_k\in N_{\Omega_1}(x_k)\cap (-N_{\Omega_2}(x_k))$ such that $\xi_k\neq 0$. By the fact that $x_k\in \Omega_1\cap \Omega_2\cap V_{k,\frac{1}{k}}(\infty; u)$, we have 
	$$\left\|\frac{x_k}{\|x_k\|}-u\right\|<\frac{1}{k}.$$
	Hence, $x_k\in \Omega_1\cap \Omega_2$ and $\tfrac{x_k}{\|x_k\|}\to u$ as $k\to\infty$. By \cite[Proposition 2.1.9]{AT03}, one has $$u\in (\Omega_1\cap \Omega_2)^\infty\subset \Omega_1^\infty\cap \Omega_2^\infty.$$ 
	Let $\eta_k:=\frac{\xi_k}{\|\xi_k\|}$. Then,  $\eta_k\in N_{\Omega_1}(x_k)\cap (-N_{\Omega_2}(x_k))$ and $\|\eta_k\|=1$. Without loss of generality, we may assume that $\eta_k\to \eta$ with $\|\eta\|=1$.  By  Proposition \ref{pro1}, we obtain that
	$$\eta\in N_{\Omega_1}(\infty;u)\cap (-N_{\Omega_2}(\infty;u)),$$
	contrary to \eqref{37a}.
	
	To prove \eqref{37b}, we let any $\xi \in N_{\Omega_1\cap \Omega_2}(\infty;u)$. From Proposition \ref{pro1} it follows that there exist sequences $x_k\in \Omega_1\cap \Omega_2$ and $\xi_k\in N_{\Omega_1\cap \Omega_2}(x_k)$ such that $x_k\xrightarrow{u}~\infty$ and $\xi_k\to~\xi$. By the fact $x_k\xrightarrow{u}\infty$,  assume that $x_k\in V_{R,\frac{1}{R}}(\infty; u)$  for every $k$ large enough. By using \eqref{37c}, we obtain that 
	$$N_{\Omega_1}(x_k)\cap (-N_{\Omega_2}(x_k))=\{0\}.$$ 
	By \cite[Theorem 2.16]{Mordukhovich2018}, we have
	$$N_{\Omega_1\cap \Omega_2}(x_k)\subset N_{\Omega_1}(x_k)+N_{\Omega_2}(x_k).$$
	Then, there exist $\eta_k\in N_{\Omega_1}(x_k)$ and $\zeta_k\in N_{\Omega_2}(x_k)$ satisfying $\xi_k=\eta_k+\zeta_k$. 
	
	Suppose that the sequence $\eta_k$ is unbounded. Hence, the sequence $\zeta_k$ is also unbounded. Without  loss of generality, we may assume that $\|\eta_k\|\to \infty$ and $\|\zeta_k\|\to\infty$.  By the fact that $\xi_k=\eta_k+\zeta_k\to \xi$, we obtain that $\frac{\eta_k}{\|\zeta_k\|}$ is bounded. By passing to subsequences if necessary, we assume that $\frac{\eta_k}{\|\zeta_k\|}\to  \eta$, $\frac{\zeta_k}{\|\zeta_k\|}\to \zeta$ for some $\eta$  and  $\zeta\in \R^n$ with $\|\zeta\|=1$. Dividing both sides of the equality $\xi_k=\eta_k+\zeta_k$ by $\|\zeta_k\|$, letting $k\to \infty$, and using Proposition \ref{pro1} yield 
	$$0=\eta+\zeta\in N_{\Omega_1}(\infty;u)\cap (-N_{\Omega_2}(\infty;u))$$ 
	with $\|\zeta\|=1$, contrary to \eqref{37a}. This implies that $\eta_k$ is bounded, and so is $\zeta_k$.
	
	Without loss of generality, assume that $\eta_k\to \eta$ and $\zeta_k\to \zeta$ for some $\eta, \zeta\in \R^n$. 
	Then, $\eta\in N_{\Omega_1}(\infty;u)$ and $\zeta\in N_{\Omega_2}(\infty;u)$ by Proposition \ref{pro1}. Hence, $\xi=\eta+\zeta\in N_{\Omega_1}(\infty;u)+ N_{\Omega_2}(\infty;u)$. \qed
 
\medskip
Let $\Omega_1$ and $\Omega_2$ be locally closed and unbounded subsets in $\R^n$ and $\R^m$, respectively. The  authors in \cite{Kim-Tung-Son-23} proved that $$N_{\Omega_1}(\infty) \times N_{\Omega_2}(\infty) \subset N_{\Omega_1\times \Omega_2}(\infty).$$ 
However, this relation does not hold for directional normal cones at infinity. To illustrate, let us consider the following examples.
\begin{example}
	$(a)$ Let $\Omega_1=\Omega_2=\R_+$. We can check that 
	$$N_{\Omega_1}(\infty;1)=\{0\}, \; \; N_{\Omega_2}(\infty;0)=\emptyset,$$
	and 
	$$N_{\Omega_1\times \Omega_2}(\infty;(1,0))=\{0\}\times  \R_-.$$
	Thus $ N_{\Omega_1\times \Omega_2}(\infty; (1,0)) \nsubseteq  N_{\Omega_1}(\infty;1) \times N_{\Omega_2}(\infty;0).$
	
	$(b)$ Let $\Omega_1:=\{2^n\;|\; n\in \N\}$ and $\Omega_2:=\{3.2^n\;|\; n\in \N\}$. We first claim that $u=(1,1) \notin (\Omega_1\times \Omega_2)^\infty$. Indeed, if otherwise, then there exist sequences $(2^{n_k}, 3. 2^{m_k})\in \Omega_1\times \Omega_2$ and $t_k\to +\infty$ such that 
	$$\lim_{k\to \infty}\frac{1}{t_k}(2^{n_k}, 3. 2^{m_k})=(1,1).$$
	It implies that $\lim_{k\to \infty} 3. 2^{m_k-n_k}=1$. This is impossible for any $n_k, m_k\in \N$. Therefore, we have $u=(1,1) \notin (\Omega_1\times \Omega_2)^\infty$ and 
	$$N_{\Omega_1\times \Omega_2}(\infty;(1,1))=\emptyset.$$
	It is easy to check that
	$N_{\Omega_1}(\infty;1)=\R$ and $N_{\Omega_2}(\infty;1)=\R$ and thus $$N_{\Omega_1}(\infty;1) \times N_{\Omega_2}(\infty;1) \nsubseteq N_{\Omega_1\times \Omega_2}(\infty;(1,1)).$$
\end{example}

\subsection{Directional Subdifferential at Infinity}
Let $f\colon\R^n\to \overline{\R}$ be an l.s.c. function and $u\in\mathbb{S}$. {\em Hereafter}, we assume that $\dom f$ is unbounded in direction $u$, i.e., $u\in(\dom f)^\infty$.   The directional subdifferentials at infinity of $f$ are introduced through the following definition. 
\begin{definition}\rm The {\em limiting} and the {\em singular  subdifferentials of $f$ in direction $u$ at infinity} are defined, respectively, by
	\begin{align*}
		\partial f(\infty; u)&:=\bigg\{\xi\in\R^n\,|\, (\xi, -1)\in \Limsup_{x\xrightarrow{u}\infty,\, r\geq f(x)}\widehat{N}_{\epi f}(x, r)\bigg\},
		\\
		\partial^\infty f(\infty; u)&:=\bigg\{\xi\in\R^n\,|\, (\xi, 0)\in \Limsup_{x\xrightarrow{u}\infty,\, r\geq f(x)}\widehat{N}_{\epi f}(x, r)\bigg\}.
	\end{align*}
\end{definition}

Based on the definition, the next proposition can be established.
\begin{proposition}\label{pro-2} For each $u\in\mathbb{S}$, one has 
	\begin{align*}  
		\partial f(\infty; u)&=\{\xi\,|\, (\xi, -1)\in \mathcal{N}_u\},
		\\
		\partial^\infty f(\infty; u)&=\{\xi\,|\, (\xi, 0)\in \mathcal{N}_u\},
	\end{align*}
	where $\mathcal{N}_u:=\displaystyle\Limsup_{x\xrightarrow{u}\infty} N_{\epi f}(x, f(x))$.
\end{proposition}
\noindent{\it Proof\,} By definition, it suffices to prove that
	\begin{equation}\label{equa-n-2}
		\displaystyle\Limsup_{x\xrightarrow{u}\infty} N_{\epi f}(x, f(x))=\Limsup_{x\xrightarrow{u}\infty,\, r\geq f(x)}\widehat{N}_{\epi f}(x, r).
	\end{equation}
	Indeed, analysis similar to that in the proof of Proposition \ref{pro1} shows that
	$$\Limsup_{x\xrightarrow{u}\infty,\, r\geq f(x)}\widehat{N}_{\mathrm{epi}\, f}(x, r)=\Limsup_{x\xrightarrow{u}\infty,\, r\geq f(x)}{N}_{\epi f}(x, r)\supset \displaystyle\Limsup_{x\xrightarrow{u}\infty} N_{\epi f}(x, f(x)).$$
	On the other hand, for any $x\in \dom f$ and $r\geq f(x)$, we see that
	\begin{equation}\label{equa-n-1}
		T_{\epi f}(x, f(x))\subset T_{\epi f}(x, r).
	\end{equation}
	Indeed, let any $v\in T_{\epi f}(x, f(x))$. Then, there exist sequences $t_k\downarrow 0$ and $(x_k, r_k)\to (x, f(x))$ such that $(x_k, r_k)\in\epi f$ for all $k\in\N$,   $\tfrac{(x_k, r_k)-(x, f(x))}{t_k}\to v$ as $k\to\infty$. Put $r^\prime_k:=r+(r_k-f(x))$. Then $r^\prime_k\to r$ and
	\begin{equation*}
		f(x_k)\leq r_k\leq r_k+(r-f(x))=r^\prime_k.
	\end{equation*}
	Hence, $(x_k, r_k^\prime)\in\epi f$ for all $k\in\N$ and 
	\begin{equation*}
		\frac{(x_k, r^\prime_k)-(x, r)}{t_k}=\frac{(x_k, r_k)-(x, f(x))}{t_k}\to v
	\end{equation*}
	as $k\to\infty$. This means that $v\in T_{\epi f}(x, r)$, as required. By \eqref{equa-n-1}, one has
	\begin{equation*}
		\widehat{N}_{\epi f}(x, r)\subset {N}_{\epi f}(x, f(x)). 
	\end{equation*}
	Hence
	$$\displaystyle\Limsup_{x\xrightarrow{u}\infty} N_{\epi f}(x, f(x))\supset \Limsup_{x\xrightarrow{u}\infty,\, r\geq f(x)}\widehat{N}_{\epi f}(x, r)$$
	and so \eqref{equa-n-2} is true. The proof is complete.\qed

\medskip 
The subsequent result is of central importance for the forthcoming analysis.
\begin{proposition}\label{pro-3}  For each $u\in\mathbb{S}$, we have 
	\begin{align}
		\partial f(\infty; u)&=\Limsup_{x\xrightarrow{u}\infty}{\partial}f(x)=\Limsup_{x\xrightarrow{u}\infty}\widehat{\partial}f(x), \label{equ-n-3}
		\\
		\partial^\infty f(\infty; u)&=\Limsup_{x\xrightarrow{u}\infty, r\downarrow 0}r{\partial}f(x)\supset\Limsup_{x\xrightarrow{u}\infty}{\partial}^\infty f(x). \label{equ-n-4}
	\end{align}
\end{proposition}

\noindent{\it Proof\,} To prove \eqref{equ-n-3}, we need to show that 
	\begin{equation*}
		\Limsup_{x\xrightarrow{u}\infty}\widehat{\partial}f(x)\subset  \Limsup_{x\xrightarrow{u}\infty}{\partial}f(x) \subset \partial f(\infty; u) \subset  \Limsup_{x\xrightarrow{u}\infty}\widehat{\partial}f(x). 
	\end{equation*}
	The first inclusion follows directly from the fact that $\widehat{\partial}f(x)\subset {\partial}f(x)$ for all $x\in~\dom f$. To prove the second inclusion, take any $\xi \in \Limsup_{x\xrightarrow{u}\infty}{\partial}f(x)$. Then there exist sequences $x_k\xrightarrow{u}\infty$ and $\xi_k\in \partial f(x_k)$ with $\xi_k\to \xi$ as $k\to\infty$. Hence $(\xi_k, -1)\in N_{\epi f}(x_k, f(x_k))$. Clearly, $(\xi_k, -1)\to (\xi, -1)$. Thus, $(\xi, -1)\in \mathcal{N}_u$, or, equivalently, $\xi\in \partial f(\infty; u)$. For the third inclusion,  let $\xi\in \partial f(\infty; u)$. Then, there exist   $x_k\xrightarrow{u}\infty$ and $(\xi_k, -r_k)\to (\xi, -1)$  with $(\xi_k, -r_k)\in N_{\epi f}(x_k, f(x_k))$ for all $k\in\N$. By definition of the limiting normal cone, for each $k\in\N$, there exist $z_k\in\R^n$ and $(\eta_k, -\alpha_k)\in\widehat{N}_{\epi f}(z_k, f(z_k))$ such that
	\begin{equation*}
		\|z_k-x_k\|\leq \frac{1}{k}, \|(\eta_k, -\alpha_k)- (\xi_k, -r_k)\|\leq \frac{1}{k}.
	\end{equation*}
	Using an argument similar as in the proof of Proposition \ref{pro1}, one has $z_k\xrightarrow{u}\infty$ and $(\eta_k, -\alpha_k)\to (\xi, -1)$. Hence, $\alpha_k>0$ for all $k$ large enough. This implies that $(\tfrac{\eta_k}{\alpha_k}, -1)\in \widehat{N}_{\epi f}(z_k, f(z_k))$ and so $\tfrac{\eta_k}{\alpha_k}\in \widehat{\partial} f(z_k)$ for all $k$ large enough. Clearly, $\tfrac{\eta_k}{\alpha_k}\to \xi$ as $k\to\infty$. Hence, $\xi \in \Limsup_{x\xrightarrow{u}\infty}\widehat{\partial}f(x)$, as required.
	
	We are now in the position to prove the first equality in \eqref{equ-n-4}. Take any $\xi\in \partial^\infty f(\infty; u)$. Then, by Proposition \ref{pro-2}, there exist sequences $x_k\xrightarrow{u}\infty$ and $(\xi_k, -r_k)\in N_{\epi f}(x_k, f(x_k))$ with $(\xi_k, -r_k)\to (\xi, 0)$.  By \cite[Proposition 1.17]{Mordukhovich2018}, one has $r_k\geq 0$ for all $k\in\N$. If $r_k=0$ for finitely many $k\in \N$, then $r_k>0$ for all $k$ large enough. Then, $(\tfrac{1}{r_k}\xi_k, -1)\in N_{\epi f}(x_k, f(x_k))$  for all $k$ large enough and 
	\begin{equation*}
		\xi=\lim_{k\to\infty}\xi_k=\lim_{k\to\infty}r_k\big(\frac{1}{r_k}\xi_k\big)\in \Limsup_{x\xrightarrow{u}\infty, r\downarrow 0}r{\partial}f(x).
	\end{equation*}
	Now, we consider the case where $r_k=0$ for infinitely many $k\in\N$. By passing to subsequences if necessary, we may assume that $r_k=0$ for all $k\in\N$. Hence, $\xi_k\in \partial^\infty f(x_k)$. By definition,  for each $k\in\N$ there exist $z_k\in\R^n$, $\eta_k\in\widehat{\partial}f(z_k)\subset \partial f(z_k)$, and $\alpha_k\in (0, \tfrac{1}{k})$ such that
	\begin{equation*}
		\|z_k-x_k\|<\frac{1}{k}, \|\alpha_k\eta_k-\xi_k\|<\frac{1}{k}.
	\end{equation*}
	This and the fact that $x_k\xrightarrow{u}\infty$ imply that $z_k\xrightarrow{u}\infty$, $\alpha_k\downarrow 0$, and $\alpha_k\eta_k\to \xi$.   Thus $\xi\in \Limsup_{x\xrightarrow{u}\infty, r\downarrow 0}r{\partial}f(x)$, as required. 
	
	Conversely, let $\xi\in \Limsup_{x\xrightarrow{u}\infty, r\downarrow 0}r{\partial}f(x)$. Then there exist $x_k\xrightarrow{u}\infty$, $r_k\downarrow 0$, and $\xi_k\in \partial f(x_k)$ such that $r_k\xi_k\to \xi$ as $k\to\infty$. This implies that $(\xi_k, -1)\in N_{\epi f}(x_k, f(x_k))$ and hence, $(r_k\xi_k, -r_k)\in N_{\epi f}(x_k, f(x_k))$ and  $(r_k\xi_k, -r_k)\to (\xi, 0)$ with $x_k\xrightarrow{u} \infty$. This means that $(\xi, 0)\in\mathcal{N}_u$ and we therefore get $\xi\in \partial^\infty f(\infty; u)$ due to Proposition \ref{pro-2}.
	
	To finish the proof, take any $\xi\in \Limsup_{x\xrightarrow{u}\infty}{\partial}^\infty f(x)$. Then there exist $x_k\xrightarrow{u}\infty$ and $\xi_k\in {\partial}^\infty f(x_k)$ such that $\xi_k\to \xi$. Hence, for each $k\in\N$ there are $z_k\in\R^n$, $\eta_k\in\widehat{\partial}f(z_k)\subset \partial f(z_k)$, and $\alpha_k\in(0,\tfrac{1}{k})$ such that
	\begin{equation*}
		\|z_k-x_k\|<\frac{1}{k},\ \  \|\alpha_k\eta_k-\xi_k\|<\frac{1}{k}.
	\end{equation*}
	This implies that $z_k\xrightarrow{u}\infty$, $\alpha_k\downarrow 0$, and $\alpha_k\eta_k\to \xi$.  Thus we obtain that $\xi\in \Limsup_{x\xrightarrow{u}\infty, r\downarrow 0}r{\partial}f(x)$. The proof is complete. \qed
 
\medskip
\begin{remark}\rm 
	For an unbounded subset $\Omega\subset\R^n$, we have
	\begin{equation*}
		\partial\delta_\Omega(\infty; u)=\partial^\infty\delta_\Omega(\infty; u)=N_\Omega(\infty; u)\ \ \forall u\in\Omega^\infty\cap\mathbb{S}.
	\end{equation*}
\end{remark}

The following result is a direct consequence of Propositions \ref{pro-3} and \ref{Proposition-3.1}. 
\begin{proposition}\label{Pro-n-2} The following assertions hold:
	\begin{enumerate}[\rm(i)]
		\item $\partial f(\infty)=\cup_{u\in \mathbb{S}}\,\partial f(\infty; u)$.
		\item $\partial^\infty f(\infty)=\cup_{u\in \mathbb{S}}\,\partial^\infty f(\infty; u)$.
	\end{enumerate}
	
\end{proposition}

The next result establishes the nonemptiness of directional subdifferentials at infinity.
\begin{proposition} \label{pro48} $\partial f(\infty; u)\cup (\partial^\infty f(\infty; u)\setminus\{0\})\neq \emptyset$.
\end{proposition}
\noindent{\it Proof\,} 	We first show that $$\mathcal{N}_u=\displaystyle\Limsup_{x\xrightarrow{u}\infty} N_{\epi f}(x, f(x))\neq \{0\}.$$
	Indeed, by the assumption that $f$ is l.s.c., we obtain that $\epi f$
	is closed.  Since $\dom f$ is unbounded in direction $u$, there exists $x_k\in~\dom f$ such that $x_k\xrightarrow{u}\infty$. For each $k\in\N$, by \cite[Proposition 1.2]{Mordukhovich2018},  there exists $\xi_k\in N_{\epi f}(x_k, f(x_k))$ with $\|\xi_k\|=1$. Then, there is  a subsequence $\xi_{k_m}$ of $\xi_k$ such that $\xi_{k_m}\to \bar \xi$ for some $\bar \xi\in \R^n$ with $\|\bar \xi\|=1$. This implies that $\bar\xi\in \mathcal{N}_u$.
	
	The conclusion follows from the definitions of the limiting and the singular subdifferentials at infinity in the direction and Proposition \ref{pro-2}. \qed

\medskip
The calculus of the limiting and singular subdifferentials at infinity in a given direction for the sum of two functions is established as follows.
\begin{proposition} \label{pro49}
	Let $f_1, f_2: \R^n \to \bar \R$ be l.s.c. functions and the following condition is satisfied
	\begin{equation}\label{d5}
		\partial^\infty f_1(\infty; u)\cap (-\partial^\infty f_2(\infty; u))=\{0\}.
	\end{equation}
	Then, one has
	\begin{equation}\label{5a}
		\partial (f_1+f_2)(\infty; u)\subset \partial f_1(\infty; u)+\partial f_2(\infty; u)
	\end{equation}
	and 
	\begin{equation}\label{5b}
		\partial^\infty (f_1+f_2)(\infty; u)\subset \partial^\infty f_1(\infty; u)+\partial^\infty f_2(\infty; u).
	\end{equation}
\end{proposition}
\noindent{\it Proof\,} We first prove  that there exists a constant  $R>0$  satisfying
	\begin{equation}\label{d6}
		\partial^\infty f_1(x)\cap (-\partial^\infty f_2(x))=\{0\} \ \ \forall x\in V_{R,\frac{1}{R}}(\infty; u).
	\end{equation}
	Indeed, on the contrary, suppose that for each $k\in\N$ there  exist $x_k\in V_{k,\frac{1}{k}}(\infty; u)$ and $\xi_k\in \partial^\infty f_1(x_k)\cap (-\partial^\infty f_2(x_k))$ such that $\xi_k\neq 0$. Since $x_k\in V_{k,\frac{1}{k}}(\infty; u)$, we have 
	$$\left\|\frac{x_k}{\|x_k\|}-u\right\|<\frac{1}{k}.$$
	Hence, one has $\tfrac{x_k}{\|x_k\|}\to u$ as $k\to\infty$. Put $\eta_k:=\frac{\xi_k}{\|\xi_k\|}$. Then, $\|\eta_k\|=1$ and  $\eta_k\in \partial^\infty f_1(x_k)\cap (-\partial^\infty f_2(x_k))$. Without loss of generality, we may assume that $\eta_k\to \eta$ with $\|\eta\|=1$. According to Proposition \ref{pro-3}, we have
	$$\eta\in \partial^\infty f_1(\infty; u)\cap (-\partial^\infty f_2(\infty; u)),$$
	contrary to \eqref{d5}.
	
	Let any $\xi\in \partial (f_1+f_2)(\infty; u)$. By Proposition \ref{pro-3}, there exist $x_k\in \R^n$ and  $\xi_k\in \partial^\infty (f_1+f_2)(x_k)$ such that $\xi_k\to \xi$ and $x_k\xrightarrow{u}\infty$. By the fact $x_k\xrightarrow{u}\infty$, we can assume that $x_k\in V_{R,\frac{1}{R}}(\infty; u)$  for every $k$ large enough. Using the condition \eqref{d6}, we have  
	$$\partial^\infty f_1(x_k)\cap (-\partial^\infty f_2(x_k))=\{0\}.$$
	By applying \cite[Theorem 2.19]{Mordukhovich2018},  one has $\xi_k=\eta_k+\zeta_k$ with $\eta_k\in \partial f_1(x_k)$ and $\zeta_k\in \partial f_2(x_k)$ for all $k$ large enough.  If $\eta_k$ is unbounded, then $\zeta_k$ is also unbounded. Without  loss of generality, we may assume that $\|\eta_k\|\to \infty$ and $\|\zeta_k\|\to\infty$.  By the fact that $\xi_k=\eta_k+\zeta_k\to \xi$, we obtain that $\frac{\eta_k}{\|\zeta_k\|}$ is bounded. By passing to subsequences if necessary, we assume that $\frac{\eta_k}{\|\zeta_k\|}\to  \eta$, $\frac{\zeta_k}{\|\zeta_k\|}\to \zeta$ for some $\eta$  and  $\zeta\in \R^n$ with $\|\zeta\|=1$. Dividing both sides of the equality $\xi_k=\eta_k+\zeta_k$ by $\|\zeta_k\|$, letting $k\to \infty$, and using Proposition \ref{pro-3} produce  $0=\eta+\zeta\in \partial^\infty f_1(\infty; u)+\partial^\infty f_2(\infty; u)$ with $\|\zeta\|=1$, which contradicts \eqref{d5}. This implies that $\eta_k$ is bounded, and so is $\zeta_k$. Without loss of generality, assume that $\eta_k\to \eta$ and $\zeta_k\to \zeta$ for some $\eta, \zeta\in \R^n$. 
	Then, $\eta\in \partial f_1(\infty;u)$ and $\zeta\in \partial f_2(\infty;u)$ by Proposition \ref{pro-3}. Thus, $\xi=\eta+\zeta\in \partial f_1(\infty; u)+\partial f_2(\infty; u)$. Therefore, the inclusion \eqref{5a} is proven.
	The proof of the inclusion \eqref{5b} is similar and so is omitted. \qed

\medskip
The following result presents the calculus of the limiting and the singular subdifferentials at infinity in a given direction for the maximum of two functions.
\begin{proposition} \label{pro411}
	Let $f_1, f_2: \R^n \to \bar \R$ be l.s.c. functions and assume that   condition \eqref{d5}  is satisfied.
	Then, one has
	\begin{equation}\label{6a}
		\partial \big(\max \{f_1, f_2\}\big)(\infty; u)\subset \bigcup_{\substack{\lambda_1, \lambda_2\in [0,1],
		\\ \lambda_1+\lambda_2=1}}\big\{\lambda_1\circ\partial f_1(\infty; u)+\lambda_2\circ\partial f_2(\infty; u)\big\}
	\end{equation}
	and 
	\begin{equation}\label{6b}
		\partial^\infty \big(\max \{f_1, f_2\}\big)(\infty; u)\subset \partial^\infty f_1(\infty; u)+\partial^\infty f_2(\infty; u),
	\end{equation}
	where $\lambda\circ\partial f(\infty; u)$ is defined as $\lambda\partial f(\infty; u)$ if $\lambda>0$ and as $\partial^\infty f(\infty; u)$ if $\lambda=0$.
\end{proposition}
\noindent{\it Proof\,} 	Let any $\xi \in \partial (\max \{f_1, f_2\})(\infty; u)$. Then, by Proposition \ref{pro-3}, there exist $x_k\in \R^n$ and  $\xi_k\in \partial (\max \{f_1, f_2\})(x_k)$ such that $x_k\xrightarrow{u}\infty$ and $\xi_k\to \xi$. By a similar argument as in the proof of Proposition \ref{pro49}, we have
	\begin{equation*}
		\partial^\infty f_1(x_k)\cap (-\partial^\infty f_2(x_k))=\{0\}.
	\end{equation*}
	Combining this with \cite[Theorem 4.10]{Mordukhovich2018}, we obtain that
	$$\xi_k\in \lambda_{1k}\circ\partial f_1(x_k)+\lambda_{2k}\circ\partial f_2(x_k)$$
	for some $\lambda_{1k}, \lambda_{2k}\in [0, 1]$ satisfying $\lambda_{1k}+ \lambda_{2k}=1$. Without loss the generality, we may assume that $\lambda_{1k}\to \lambda_1$ and $\lambda_{2k}\to\lambda_2$ for some $\lambda_1, \lambda_2\in [0, 1]$ satisfying $\lambda_1+ \lambda_2=1$. We consider the following cases:
	
	{\it Case 1: $\lambda_1=0$ and $\lambda_2=1$}. Put $K:=\{k: \lambda_{1k}=0\}$. There are two subcases to be considered.
	
	{\it Case 1.1: $K$ is finite.} Then, for $k$ large enough, there exist $\eta_k\in \partial f_1(x_k)$ and $\zeta_k\in \partial f_2(x_k)$ satisfying $\xi_k=\lambda_{1k}\eta_k+\lambda_{2k}\zeta_k$. By a same argument as in the proof of Proposition \ref{pro49}, the sequence $\zeta_k$ is bounded. We can assume that $\zeta_k\to\zeta$ with $\zeta\in \partial f_2(\infty;u)$. Then, one has $\lambda_{1k}\eta_k\to \xi-\zeta$. Since $\lambda_{1k}\to 0$ and $\eta_k\in \partial f_1(x_k)$, by Proposition \ref{pro-3}, we obtain that $\xi-\zeta\in \partial^\infty f_1(\infty;u)$. Hence, 
	$$\xi\in \partial^\infty f_1(\infty;u)+\partial f_2(\infty;u).$$ 
	
	{\it Case 1.2: $K$ is infinite}. Without loss of generality,  assume that $\lambda_{1k}=0$ for all $k\in\N$ and $\xi_k=\tilde{\eta}_k+\lambda_{2k}\zeta_k$ with $\tilde{\eta}_k\in \partial^\infty f_1(x_k)$ and $\zeta_k\in \partial f_2(x_k)$. Using the same argument as in the proof of Proposition \ref{pro49} and passing subsequences if necessary, we can assume that $\tilde{\eta}_k\to \tilde{\eta}$ and $\zeta_k\to \zeta$ for some $\tilde{\eta}\in \partial^\infty f_1(\infty;u)$ and $\zeta\in \partial f_2(\infty;u)$. This implies that 
	$$\xi \in  \partial^\infty f_1(\infty;u)+\partial f_2(\infty;u).$$ 
	
	{\it Case 2: $\lambda_1=1$ and $\lambda_2=0$}. The proof follows from the same argument as in Case 1.
	
	{\it Case 3: $\lambda_1, \lambda_2\in (0, 1)$.} Then,  assume that $\lambda_{1k}>0$ and $\lambda_{2k}>0$ for all $k$ and $\xi_k=\lambda_{1k}\eta_k+\lambda_{2k}\zeta_k$ with $\eta_k\in \partial f_1(x_k)$ and $\zeta_k\in \partial f_2(x_k)$. Using the same analysis as in Proposition \ref{pro49}, we can assume that $\eta_k\to \eta$ and $\zeta_k\to \zeta$ for some $\eta\in \partial f_1(\infty;u)$ and $\zeta\in \partial f_2(\infty;u)$. Then, we obtain that 
	$$\xi \in  \lambda_1\partial f_1(\infty;u)+\lambda_2\partial f_2(\infty;u).$$
	
	By the above cases, \eqref{6a} is proven. The proof of \eqref{6b} is similar. \qed 
	
\medskip
From Proposition \ref{pro-3}, the following result can be derived.
\begin{proposition} \label{pro412}
	Let $f_1, f_2: \R^n \to \bar \R$ be l.s.c. functions. 
	Then, one has
	\begin{equation}\label{8a}
		\partial (\min \{f_1, f_2\})(\infty; u)\subset\partial f_1(\infty; u)\cup\partial f_2(\infty; u)
	\end{equation}
	and 
	\begin{equation}\label{8b}
		\partial^\infty (\min \{f_1, f_2\})(\infty; u)\subset\partial^\infty f_1(\infty; u)\cup\partial^\infty f_2(\infty; u).
	\end{equation}
\end{proposition}
\noindent{\it Proof\,} Let any $\xi \in \partial (\min \{f_1, f_2\})(\infty; u)$. According to Proposition \ref{pro-3}, there exist $x_k\in~\R^n$ and  $\xi_k\in \partial (\min \{f_1, f_2\})(x_k)$ such that $x_k\xrightarrow{u}\infty$ and $\xi_k\to \xi$. According to \cite[Proposition 4.9]{Mordukhovich2018}, we have $\xi_k\in \partial f_1(x_k)\cup \partial f_2(x_k)$. Using Proposition \ref{pro-3}, we obtain that $\xi \in \partial f_1(\infty; u)\cup\partial f_2(\infty; u)$ and the inclusion \eqref{8a} is proven. 
	
The proof of \eqref{8b} is analogous to that of \eqref{8a}, and is therefore omitted.\qed 

\medskip
We now present formulae to estimate the directional partial subdifferentials at infinity of two variables functions.   Let $u\in \mathbb{S}$ and  consider a proper l.s.c. function $F:\R^n\times \R^m \to \overline{\R}, \, (x,y)\mapsto F(x,y)$. Fix $\bar y\in \R^m$ and denote, respectively, by $\partial F_x(\infty;  \bar y, u)$ and $\partial^\infty F_x(\infty; \bar y, u)$ the limiting and the  singular  subdifferentials of $F(\cdot, \bar y):\R^n\to \overline{\R}, \, x\mapsto F(x,\bar y)$ at infinity in   direction $u$.

The following  establishes estimates for $\partial F_x(\infty; \bar y, u)$ and $\partial^\infty F_x(\infty; \bar y, u)$.
\begin{proposition} \label{pro413} 
	If the following condition is satisfied
	\begin{equation}\label{9a}
		(0, \eta) \in \partial^\infty F(\infty;  (u, 0)) \; \Rightarrow \; \eta =0,
	\end{equation}
	then 
	\begin{equation}\label{9b}
		\partial F_x(\infty; \bar y, u) \subset \{\xi \in \R^n \;|\; \exists \eta\in \R^m \;  with \; (\xi,\eta)\in \partial F(\infty; (u,0))\}
	\end{equation}
	and 
	\begin{equation}\label{9c}
		\partial^\infty F_x(\infty;  \bar y, u) \subset \{\xi \in \R^n \;|\; \exists \eta\in \R^m \;  with \; (\xi,\eta)\in \partial^\infty F(\infty;(u,0))\}.
	\end{equation}
\end{proposition}
\noindent{\it Proof\,} We first claim  that there exists a constant  $R>0$  such that for every $x\in V_{R,\frac{1}{R}}(\infty; u)$ satisfying
	\begin{equation}\label{9d}
		(0, \eta) \in \partial^\infty F(x, \bar y) \; \Rightarrow \; \eta =0.
	\end{equation}
	Indeed, on the contrary, suppose that for each $k\in\N$ there  exist $x_k\in V_{k,\frac{1}{k}}(\infty; u)$ and $(0, \eta_k) \in \partial^\infty F(x_k, \bar y)$ with $\eta_k\neq 0$. By the fact that $x_k\in V_{k,\frac{1}{k}}(\infty; u)$, we have 
	$$\left\|\frac{x_k}{\|x_k\|}-u\right\|\leq \frac{1}{k}.$$
	This implies $x_k\xrightarrow{u}\infty$. Thus 
	$$\left\|\frac{(x_k,\bar y)}{\|(x_k,\bar y)\|}-(u,0)\right\|\leq\left\|\frac{x_k}{\|x_k\|+\|\bar y\|}-u\right\|+\left\|\frac{\bar y}{\|x_k\|+\|\bar y\|}\right\|\to 0$$
	as $k\to \infty$. Hence, $(x_k,\bar y)\xrightarrow{(u, 0)}\infty$ as $k\to\infty$. Put $\zeta_k:=\frac{\eta_k}{\|\eta_k\|}$. Then, $\|\zeta_k\|=1$ and $(0, \zeta_k) \in \partial^\infty F(x_k, \bar y)$. Without loss of generality, we may assume that $\zeta_k\to \zeta$ with $\|\zeta\|=1$. By using Proposition \ref{pro-3}, one has
	$$(0, \zeta) \in \partial^\infty F(\infty;(u,0)),$$
	contrary to \eqref{9a}.

	To prove the inclusion \eqref{9b}, we  take any $\xi \in \partial F_x(\infty; \bar y, u)$. By Proposition \ref{pro-3}, there are sequences $x_k\xrightarrow{u}\infty$  and $\xi_k\in \partial F_x(x_k, \bar y)$
	such that $\xi_k\to \xi$. For all $k$  large enough, the condition \eqref{9d} holds at $x=x_k$. By \cite[Corollary 10.11]{Rockafellar1998}, we obtain  that
	$$\partial F_x(x_k, \bar y) \subset \{\xi \in \R^n |\ \exists \eta\in \R^m \ \  \text{with} \ \ (\xi,\eta)\in \partial F(x_k,\bar y)\}.$$
	Then, there is $\eta_k\in \R^m$ such that $(\xi_k,\eta_k)\in \partial F(x_k,\bar y)$. 
	
	We claim that the sequence $\eta_k$ is bounded. Indeed, if otherwise, we can assume that $\eta_k\to \infty$  and $\frac{\eta_k}{\|\eta_k\|}\to \eta$ with $\|\eta\|=1$. Then, one has
	$$\bigg(\frac{\xi_k}{\|\eta_k\|}, \frac{\eta_k}{\|\eta_k\|}\bigg)\in \frac{1}{\|\eta_k\|} \partial F(x_k,\bar y).$$
	By Proposition \ref{pro-3}, we obtain that  $(0, \eta) \in \partial^\infty F(\infty;(u,0))$, contrary to \eqref{9a}. Thus,  the sequence $\eta_k$ is bounded.
	
	Finally, we assume that $\eta_k\to\eta$. As $x_k\xrightarrow{u}\infty$, we get $(x_k,\bar y)\xrightarrow{(u,0)}\infty$ and $(\xi_k,\eta_k)\to (\xi, \eta)$. By Proposition \ref{pro-3}, one has  $(\xi,\eta)\in \partial F(\infty; (u,0))$ and the inclusion \eqref{9b} is proven. 
	
	The proof of \eqref{9c} is  similar to that of \eqref{9b} and is therefore omitted. \qed
 
\subsection{Directional Lipschitzness at Infinity}
This subsection provides a characterization of directional Lipschitzness at infinity via directional subdifferentials at infinity.
\begin{definition} \rm Let $f\colon\R^n\to \R$ be a real-valued function and $u\in\R^n$. We say that $f$ is {\em Lipschitz at infinity in  direction $u$} if there exist $L>0$, $R>0$, and $\delta>0$ such that
	\begin{equation*}
		|f(x)-f(y)|\leq L\|x-y\|\ \ \forall x,y\in V_{R,\delta}(\infty; u).
	\end{equation*}
\end{definition}

A necessary and sufficient condition for  directional Lipschitzness at infinity is given below.
\begin{proposition}\label{Pro-n-1} The function $f$ is Lipschitz at infinity in  direction $u$ if and only if $\partial^\infty f(\infty; u)=\{0\}$. In this case, the set  $\partial f(\infty; u)$ is nonempty and compact.
\end{proposition}
\noindent{\it Proof\,} Assume that $f$ is  Lipschitz at infinity in   direction $u$. Then there exist $L>0$, $R>0$, and $\delta>0$ such that
	\begin{equation*}
		|f(x)-f(y)|\leq L\|x-y\|\ \ \forall x,y\in V_{R,\delta}(\infty; u).
	\end{equation*}  
	Let $\xi\in \partial^\infty f(\infty; u)$. By Proposition \ref{pro-3}, there exist $x_k\xrightarrow{u}\infty$, $\xi_k\in \partial f(x_k)$, and $r_k\downarrow 0$ such that $r_k\xi_k\to \xi$ as $k\to\infty$. Since $x_k\xrightarrow{u}\infty$, there exists $k_0\in\N$ such that $x_k\in V_{R, \tfrac{\delta}{2}}(\infty; u)$ for all $k\geq k_0$. This means that $x_k\in\mathrm{int}\, V_{R, \delta}(\infty; u)$ and so $f$ is locally Lipschitz around $x_k$ with constant $L$ for all $k\geq k_0$. By \cite[Theorem 1.22]{Mordukhovich2018}, $\|\xi_k\|\leq L$ for all $k\geq k_0$. This implies that $r_k\xi_k\to 0$ as $k\to\infty$ and so $\xi=0$. Hence, $\partial^\infty f(\infty; u)=\{0\}$.
	
	We now assume that $\partial^\infty f(\infty; u)=\{0\}$. Then there exist $L>0$, $R>0$, and $\delta>0$ such that for all $x\in V_{R,\delta}(\infty; u)$ and $\xi\in\partial f(x)$, $\|\xi\|\leq L$. Indeed, if otherwise, then for each $k\in\mathbb{N}$, there exist $x_k\in V_{k, \tfrac{1}{k}}(\infty; u)$ and $\xi_k\in\partial f(x_k)$ such that $\|\xi_k\|>k$. Hence, $x_k\xrightarrow{u}\infty$, $\|\xi_k\|\to\infty$, and $\tfrac{\xi_k}{\|\xi_k\|}\to \xi$ with $\|\xi\|=1$ (by passing to a subsequence  if necessary). By Proposition \ref{pro-3}, $\xi\in\partial^\infty f(\infty; u)\setminus\{0\}$, a contradiction.
	
	Now, by \cite[Theorem I]{Whitney}, there exists a $C^\infty$-function $\varphi\colon\R^n\to [0, 1]$ satisfying
	\begin{equation*}
		\varphi(x)=
		\begin{cases}
			1, \ \ &\text{if}\ \ \|x\|\geq 3R,
			\\
			0, \ \ &\text{if}\ \ \|x\|\leq 2R.
		\end{cases}
	\end{equation*}
	Let $\bar f$ be the function defined by $\bar f(x)=f(x)\varphi(x)$ for all $x\in\R^n$. It is easy to check that $\bar f$ is locally Lipschitz on 
	$$V_{3R,\delta}(u):=\Big\{z\in\mathbb{B}_{3R}\,|\, \big\|z\|u\| -u\|z\|\big\|\leq \delta\|z\|\|u\|\Big\}.$$
	Hence, $\bar f$ is globally Lipschitz on $V_{3R,\delta}(u)$. By increasing $L$ if necessary, assume that $\bar f$ is locally Lipschitz with constant $L$ on $V_\delta(u):=V_{3R,\delta}(u)\cup V_{3R, \delta}(\infty; u)$. Hence, by \cite[Theorem 1.22]{Mordukhovich2018}, $\|\xi\|\leq L$ for all $\xi\in\partial\bar f(x)$ and all $x\in V_{\tfrac{\delta}{2}}(u)$. By the mean value theorem \cite[Corollary 4.14]{Mordukhovich2018}, one has
	\begin{equation*}
		|\bar f(x)-\bar f(y)|\leq \|x-y\|\sup\{\|\xi\|\,|\, \xi\in\partial\bar f(u), u\in [x, y]\}\leq L\|x-y\|
	\end{equation*}
	for all $x, y\in V_{\tfrac{\delta}{2}}(u)$. Hence,
	\begin{equation*}
		|f(x)-f(y)|=|\bar f(x)-\bar f(y)|\leq L\|x-y\|\ \ \forall x, y\in V_{3R,\tfrac{\delta}{2}}(u)
	\end{equation*}
	and so $f$ is Lipschitz at infinity in  direction $u$.
	
	We now prove that if $f$ is Lipschitz at infinity in  direction $u$, then $\partial f(\infty; u)$ is nonempty and compact. Indeed, the nonemptyness of $\partial f(\infty; u)$ follows directly from Proposition \ref{pro48} and the fact that $\partial^\infty f(\infty; u)=\{0\}$. The closedness of $\partial f(\infty; u)$ is obvious from the definition. Let any $\xi\in\partial f(\infty; u)$. Then by definition, there exist $x_k\xrightarrow{u}\infty$ and $\xi_k\in\partial f(x_k)$ with $\xi_k\to \xi$. Since $x_k\xrightarrow{u}\infty$, $x_k\in V_{3R,\tfrac{\delta}{2}}(u)$ for all $k$ large enough. By the Lipschitzness at infinity in  direction $u$ with constant $L$ of $f$, we have $\|\xi_k\|\leq L$ for all $k$ large enough. Hence, $\|\xi\|\leq L$ and so $\partial f(\infty; u)$  is compact.   \qed   
 
\medskip
Let $\Omega$ be a nonempty
closed set and let $u\in \mathbb{S}^n$. Then the distance function $d_\Omega(\cdot)$ is globally Lipschitz with constant $1$. Applying Proposition \ref{Pro-n-1}, we obtain that $\partial^\infty d_\Omega(\infty; u)=\{0\}$ and the set  $\partial d_\Omega(\infty; u)$ is nonempty and compact. We derive below the explicit formula for  $\partial d_\Omega(\infty; u)$.
\begin{proposition}
	It holds that
	\begin{equation*}
		\partial d_{\Omega}(\infty; u)=
		\begin{cases}
			(N_\Omega(\infty; u)\cap \mathbb{B})\cup \mathbb{E},\ \ &\text{if}\ \Omega \text{ is unbounded,}
			\\
			\{u\}, \ &\text{otherwise,}
		\end{cases}
	\end{equation*}
	where $$\mathbb{E}:=\Limsup_{x\xrightarrow{u}\infty, \ x\notin\Omega}\frac{x-\Pi_\Omega(x)}{d_\Omega(x)}.$$
\end{proposition}
\noindent{\it Proof\,} We first consider the case where $\Omega$ is bounded.
	Let any $\xi\in \partial d_\Omega(\infty;u)$.  Then, there exist $x_k\xrightarrow{u}\infty$ and $\xi_k\to \xi$ such that $\xi_k\in \partial d_\Omega (x_k)$. By the boundedness of $\Omega$, we can assume that $x_k\notin \Omega$ for every $k$. Applying \cite[Theorem 1.33]{Mordukhovich2018}, for each $k$, we have 
	$$\partial d_\Omega (x_k)=\frac{x_k-\Pi_\Omega(x_k)}{d_\Omega(x_k)}.$$
	Then, there exists $\bar x_k\in \Pi_\Omega(x_k)$ such that $\xi_k=\frac{x_k-\bar x_k}{\|x_k-\bar x_k\|}$. Since $x_k\xrightarrow{u}\infty$ and $\Omega$ is bounded, one has $\xi_k\to u$ as $k\to \infty$. This follows that $\xi=u$ and so $\partial d_\Omega(\infty; u)\subset \{u\}$.
	
	Since the distance function $d_\Omega\colon\R^n\to \R$, $x\mapsto \inf_{y\in\Omega}\|x-y\|$, is globally Lipschitz with constant $1$, $\partial d_\Omega(\infty; u)$ is nonempty and compact due to Proposition \ref{Pro-n-1}.  	Hence, $\partial d_\Omega(\infty; u)=\{u\}$.

	Now, we consider the case where $\Omega$ is unbounded. According to \cite[Theorem 1.33]{Mordukhovich2018}, we have
	\begin{equation*}
		\partial d_{\Omega}(x)=
		\begin{cases}
			N_\Omega(x)\cap \mathbb{B},\ \ &\text{if}\ x\in \Omega,
			\\
			\frac{x-\Pi_\Omega(x)}{d_\Omega(x)}, \ &\text{otherwise.}
		\end{cases}
	\end{equation*}
	By applying Proposition \ref{pro-3}, we obtain that 
	\begin{equation*}
		\partial d_{\Omega}(\infty; u)=
		(N_\Omega(\infty; u)\cap \mathbb{B})\cup \mathbb{E}.
	\end{equation*}
	The proof is complete. \qed
 
\medskip
The following result characterizes the directional normal cone at infinity for a set defined by equality and inequality constraints.
\begin{proposition} Let $\Omega$ be a unbounded subset in $\R^n$ and $u\in\Omega^\infty\cap\mathbb{S}$.
	Consider the following set 
	\begin{equation*}
		S:=\{x\in\Omega\;|\; g_i(x)\leq 0, i=1, \ldots, m, h_j(x)=0, j=1, \ldots, p\},
	\end{equation*}
	where $g_i, h_j\colon\R^n\to \overline{\R}$, $i=1, \ldots, m,$ $j=1, \ldots, p$, are Lipschitz at infinity in  direction $u$. If $S$ is unbounded and satisfies the limiting constraint qualification at infinity in  direction $u$, i.e., there do not exist $\lambda_i\geq 0,$ $i=1, \ldots, m,$ and $\mu_j\geq 0, j=1, \ldots, p,$ not all zero, such that 
	\begin{equation}\label{DLCQ}
		0\in  \sum_{i=1}^m\lambda_i\partial g_i(\infty; u)+\sum_{j=1}^p\mu_j\big[\partial h_j(\infty; u)\cup \partial (-h_j)(\infty; u)\big]+N_\Omega(\infty; u), \tag{LCQ$_u^\infty$}
	\end{equation}  
	then 
	\begin{equation*}
		N_S(\infty; u)\subset \mathrm{pos}\,\left\{\bigcup_{i=1}^m \partial g_i(\infty; u), \bigcup_{j=1}^p \big[\partial h_j(\infty; u)\cup  \partial  (-h_j)(\infty; u)\big] \right\}+N_\Omega(\infty; u).
	\end{equation*}
\end{proposition}
\noindent{\it Proof\,} Take any $\xi\in N_S(\infty; u)$. We consider only the case $\xi \neq 0$, as the conclusion is trivial otherwise. Put 
		\begin{equation*}
			C:=\{x\in\R^n\;|\; g_i(x)\leq 0, i=1, \ldots, m, h_j(x)=0, j=1, \ldots, p\},
		\end{equation*}
		then $S=C\cap\Omega$. From the condition \eqref{DLCQ}, it is clear that
		\begin{equation*}
			0\notin \mathrm{co}\,\left\{\bigcup_{i=1}^m \partial g_i(\infty; u), \bigcup_{j=1}^p \big[\partial h_j(\infty; u)\cup \partial  (-h_j)(\infty; u)\big] \right\}.
		\end{equation*}
		Thus, by an argument similar to that in the proof of \cite[Proposition 5.9]{Kim-Tung-Son-23}, we obtain
		\begin{equation}\label{equ-n-11}
			N_C(\infty; u)\subset \mathrm{pos}\,\left\{\bigcup_{i=1}^m \partial g_i(\infty; u), \bigcup_{j=1}^p \big[\partial h_j(\infty; u)\cup \partial  (-h_j)(\infty; u)\big] \right\}.
		\end{equation}
		We claim that $N_C(\infty; u)\cap [-N_\Omega(\infty; u)]=\{0\}$ and so by Proposition~\ref{pro37} and inclusion \eqref{equ-n-11}, we have
		\begin{align*}
			N_S(\infty&, u)\subset N_C(\infty; u)+N_\Omega(\infty; u)
			\\
			&\subset \mathrm{pos}\,\left\{\bigcup_{i=1}^m \partial g_i(\infty; u), \bigcup_{j=1}^p \big[\partial h_j(\infty; u)\cup  \partial  (-h_j)(\infty; u)\big] \right\}+N_\Omega(\infty; u),
		\end{align*}
		as required. Indeed, if otherwise, then there exists  $\eta\in N_C(\infty; u)\cap [-N_\Omega(\infty; u)]$ with $\eta\neq 0$. By \eqref{equ-n-11}, there are $\lambda_i\geq 0,$ $\zeta_i\in \partial g_i(\infty; u)$ $i=1, \ldots, m,$  $\mu_j\geq 0,$ $\omega_j\in [\partial h_j(\infty; u)\cup  \partial (-h_j)(\infty; u)]$, $j=1, \ldots, p,$ such that
		$$0=\sum_{i=1}^m\lambda_i\zeta_i+\sum_{j=1}^p\mu_j\omega_j-\eta.$$ 
		This and \eqref{DLCQ} imply that $\lambda_i=0$, $i=1, \ldots, m,$ and $\mu_j= 0, j=1, \ldots, p.$ Thus $\eta=0$, a contradiction. The proof is complete. \qed

\section{Applications}\label{Section-4}
In this section, employing directional normal cones at infinity together with directional limiting and singular subdifferentials at infinity, we develop several applications, including directional optimality conditions at infinity, coercivity, compactness of the global solution set, weak sharp minima at infinity, and error bounds at infinity.

\subsection{Directional Optimality Conditions at Infinity}

Let $f \colon \mathbb{R}^n \to \overline{\mathbb{R}}$ be an l.s.c. function and  $\Omega$ be a nonempty and closed subset of $\R^n$. We consider the following general optimization problem 
\begin{equation}\label{problem-0}
	\min_{x\in\Omega} f(x).\tag{P}
\end{equation}

To present a necessary optimality condition at infinity for the problem \eqref{problem-0},  we always assume that:
\\
(A$_1$) $\mathrm{dom} f \cap \Omega$ is unbounded.
\\
(A$_2$) $f$ is bounded from below on $\Omega,$ i.e., $f_* := \inf_{x\in \Omega} f(x)$ is finite. 

\medskip
The following provides a necessary optimality condition at infinity  for problem \eqref{problem-0} in a given asymptotic direction of the feasible set. 
\begin{theorem} \label{Necessary-Theorem} 
	Let $u\in\Omega^\infty\cap\mathbb{S}$.  Assume that the following condition holds:
	\begin{equation}\label{equ-n-6}
		\partial^{\infty}f(\infty;u) \cap \big(-N_{\Omega}(\infty; u) \big) =\{0\}.
	\end{equation}
	If there exists a sequence $x_k\xrightarrow{\Omega, u}\infty$ such that $f(x_k)\to f_*:=\inf_{x\in \Omega} f(x)$, 
	then
	\begin{equation*}
		0\in\partial f(\infty; u)+N_\Omega(\infty; u).
	\end{equation*}
\end{theorem}  
\noindent{\it Proof\,} 	We first consider the case that $\Omega=\mathbb{R}^n$. Then  $N_{\Omega}(\infty; u) = \{0\}.$ For each $k\in\mathbb{N}$, we have
	\begin{eqnarray*}
		f_*\leq f(x_k)\leq f_*+\left(f(x_k)-f_*+\frac{1}{k}\right). 
	\end{eqnarray*}
	Clearly, $\epsilon_k:=f(x_k)-f_*+\frac{1}{k}>0$  and $\epsilon_k\to 0$ as $k\to\infty$. Put $\lambda_k:=\sqrt{\epsilon_k}$, then by the Ekeland variational principle (Lemma~\ref{lema26}), there exists $z_k \in \mathbb{R}^n$ for $k > 0$ such that
	\begin{align*}
		&\|x_k - z_k\|   \le   \lambda_k, \\
		&f(z_k)  \le  f(x)+\lambda_k\|x-z_k\| \quad \textrm{for all } \quad x \in \mathbb{R}^n.
	\end{align*}
	The first inequality and the fact that  $x_k\xrightarrow{u}\infty$  imply that  $z_k \xrightarrow{u}\infty$. While the second inequality says that $z_k$ is a global minimizer of   $\varphi(\cdot):=f(\cdot)+\lambda_k\|\cdot-z_k\|$  on $\mathbb{R}^n$. 	By the Fermat rule (see Lemma~\ref{lema22}), we obtain
	\begin{eqnarray*}
		0 &\in& \partial \left( f(\cdot)+ \lambda_k\|\cdot-z_k\|\right)(z_k).
	\end{eqnarray*}
	By the Lipschitzness of the function $\|\cdot-z_k\|$ and the sum rule (see \cite[Theorem~2.19]{Mordukhovich2018}), we have 
	\begin{align*}
		0 &\in \partial f(z_k)+ \lambda_k\partial (\|\cdot-z_k\|)(z_k)
		\\
		&=\partial f(z_k)+ \lambda_k\mathbb{B}
	\end{align*}
	due to the fact that $\partial (\|\cdot-z_k\|)(z_k)=\mathbb{B}$. Hence, 
	$$0 \in \partial f(z_k)+ \lambda_k\mathbb{B},$$
	and so there is $\xi_k \in \partial f(z_k)$ such that $\|\xi_k\| \leq \lambda_k$. Since $\lambda_k\to 0$, by letting $k \to \infty$ and applying Proposition~\ref{pro-3}, we obtain  $ 0\in \partial f(\infty; u).$
	
	We now consider the case where $\Omega$ is an arbitrary subset of $\mathbb{R}^n$. 
	We have
	\begin{eqnarray*}
		f_*=\inf_{x \in \mathbb{R}^n} \left(f + \delta_{\Omega}\right)(x) & = & \inf_{x \in \Omega} f(x) \ > \ -\infty,
	\end{eqnarray*}
	where $\delta_{\Omega} \colon \mathbb{R}^n \to \overline{\mathbb{R}}$  stands for the indicator function of the set $\Omega.$ Clearly $f(x_k) + \delta_{\Omega}(x_k)\to f_*$. Therefore, $0 \in \partial (f + \delta_{\Omega})(\infty; u)$ (by the argument employed in the first case). This, together with Proposition~\ref{pro49} and \eqref{equ-n-6}, yields 
	$$0 \in \partial f (\infty; u) + N_{\Omega}(\infty; u).$$ 
	The proof is complete.\qed
 
\begin{remark}\rm By Propositions \ref{pro36} and \ref{Pro-n-1},  condition \eqref{equ-n-6} holds automatically when $f$ is Lipschitz at infinity in   direction $u$, or when $u\in \Omega^\infty\setminus (\bd\Omega)^\infty$.
\end{remark}

In connection with Theorem \ref{Necessary-Theorem}, the following natural question arises:

$\,$\\
$(Q)$ {\em  If condition \eqref{equ-n-6} and the following
	\begin{equation}\label{equ-n-10}
		0\notin\partial f(\infty; u)+N_\Omega(\infty; u)
	\end{equation}
	are satisfied, then there exists $\delta>0$ such that $\mathrm{argmin}\,_{x\in \Omega\cap V_{\delta}(u)} f(x)\neq \emptyset$, where
	\begin{equation*}
		V_{\delta}(u):=\Big\{x\in\mathbb{R}^n\,|\, \big\| x-u\|x\|\big\|\leq \delta\|x\|\Big\}=\mathrm{cone}\,\mathbb{B}(u, \delta)?
	\end{equation*}
}

The following example answers question $(Q)$ in the negative.
\begin{example}\label{Example-4.3}\rm Let $f\colon\mathbb{R}^2\to \mathbb{R}$ be the function defined by $f(x_1, x_2)=x_1^2+e^{x_2}$ for all $x=(x_1, x_2)\in\Omega:=\mathbb{R}^2$. It is easy to see that $\inf_{x\in\mathbb{R}^2}f(x)=0$, $\mathrm{argmin}_{x\in\mathbb{R}^2}f(x)=\emptyset$, and 
	$$0\in \partial f(\infty)=\mathbb{R}\times\{0\}.$$
	Let $u=(1, 0)$, then  an easy computation shows that $\partial f(\infty; u)=\emptyset$ and  $N_\Omega(\infty; u)=\{0\}$. Hence, conditions  \eqref{equ-n-6} and  \eqref{equ-n-10} are satisfied. For any $\delta>0$, we can check that
	\begin{equation*}
		V_\delta(u)= \mathrm{cone}\,\mathbb{B}((1, 0), \delta),\ \ \inf_{V_\delta(u)} f(x)=0,
	\end{equation*}
	and $\mathrm{argmin}\,_{x\in V_{\delta}(u)} f(x)=\emptyset$.
\end{example}

The following results provide sufficient conditions for the existence of solutions to \eqref{problem-0} over the feasible set along with a given asymptotic direction. 

\begin{corollary}  Let $u\in\Omega^\infty\cap \mathbb{S}$ and $\bar x\in\Omega\cap\dom f$. If the following conditions hold
	\begin{equation*}
		N_\Omega(\infty; u)\cap u^\perp=\{0\}, \partial^\infty f(\infty; u)\cap [-N_\Omega(\infty; u)+u^\perp]=\{0\},
	\end{equation*}  
	and 
	\begin{equation*}
		0\notin \partial f(\infty; u)+N_\Omega(\infty; u)+u^\perp,
	\end{equation*}
	then  we have $\mathrm{argmin}_{x\in \Omega_u}\, f(x)$ is nonempty, where 
	$$\Omega_u:=\Omega\cap [\bar x+\mathrm{pos}\,\{u\}].$$ 
\end{corollary}
\noindent{\it Proof\,}  Since $f$ is bounded from below on $\Omega$ and $\Omega_u\subset\Omega$, $f_*^u:=\inf_{x\in\Omega_u} f(x)$ is finite. If $\mathrm{argmin}_{x\in \Omega_u}\, f(x)$ is empty, then it follows from Theorem \ref{Necessary-Theorem} that  $0\in\partial (f+\delta_{\Omega_u})(\infty; u)$. On the other hand, since $N_{\bar x+\mathrm{pos}\,\{u\}}(\infty; u)=u^\perp$ and $N_\Omega(\infty; u)\cap u^\perp=\{0\}$, we get     
	\begin{equation*}
		N_{\Omega_u}(\infty; u)\subset N_\Omega(\infty; u)+u^\perp
	\end{equation*}
	due to Proposition \ref{pro37}. This and the fact that $$\partial^\infty f(\infty; u)\cap [-N_\Omega(\infty; u)+u^\perp]=\{0\}$$ imply that $\partial^\infty f(\infty; u)\cap [-N_{\Omega_u}(\infty; u)]=\{0\}$. Thus, by Proposition \ref{pro49}, we see that
	\begin{equation*}
		0\in\partial (f+\delta_{\Omega_u})(\infty; u)\subset \partial f(\infty; u)+N_{\Omega_u}(\infty; u) \subset \partial f(\infty; u)+N_\Omega(\infty; u)+u^\perp, 
	\end{equation*}  
	a contradiction. The proof is complete. \qed
 
\begin{corollary}\label{Cor-n-1} If $\Omega$ is convex, $u\in\Omega^\infty\cap \mathbb{S}$, and the following conditions hold
	\begin{equation}\label{equ-n-9}
		\partial^\infty f(\infty; u)\cap u^\perp=\{0\}, \partial f (\infty; u)\cap u^\perp=\emptyset,
	\end{equation}	 
	then for any $\bar x\in\Omega\cap\dom f$, the set $\mathrm{argmin}_{x\in \bar x+\mathrm{pos}\,\{u\}}\, f(x)$ is nonempty. 
\end{corollary}
\noindent{\it Proof\,}  Let any $\bar x\in\Omega\cap \dom f$. By the convexity of $\Omega$, one has $\bar x+\mathrm{pos}\,\{u\}\subset \Omega$. If  $\mathrm{argmin}_{x\in \bar x+\mathrm{pos}\,\{u\}}\, f(x)=\emptyset$, then by Theorem \ref{Necessary-Theorem}, one has
	$$0\in \partial (f+\delta_{\bar x+\mathrm{pos}\,\{u\}})(\infty; u).$$
	This and \eqref{equ-n-9} imply that $0\in \partial f(\infty; u)+u^\perp.$ Hence, $\partial f (\infty; u)\cap u^\perp$ is nonempty, a contradiction. The proof is complete. \qed
 
\begin{example}\rm Let $f$ and $\Omega$ be as in Example \ref{Example-4.3}. Then, an easy computation shows that $\partial f(\infty; u)=\emptyset$ and $\partial^\infty f(\infty; u)=\mathbb{R}\times\{0\}$  for $u=(\pm 1, 0)$. Hence, \eqref{equ-n-9} is satisfied and so  $\mathrm{argmin}_{x\in\mathrm{pos}\, u}f(x)$ is nonempty due to Corollary \ref{Cor-n-1}. In fact, it is easy to check that $\mathrm{argmin}_{x\in\mathrm{pos}\, u}f(x)=\{(0,0)\}$ for $u=(\pm 1, 0)$.
	
\end{example}

The nonemptiness and compactness of the global solution, coercivity, and weak sharp minima at infinity are investigated in the following result.
\begin{theorem}\label{Theorem3.3} Assume that  \eqref{equ-n-6}  and the following condition
	\begin{equation}\label{equ-n-5}
		0\notin\partial f(\infty; u)+N_\Omega(\infty; u)  
	\end{equation}
	hold  for all $u\in \Omega^\infty\cap \mathbb{S}$, then the following assertions hold:
	\begin{enumerate}[\rm(i)]
		\item $\Sol\eqref{problem-0}$ is nonempty and compact.
		\item Problem \eqref{problem-0} has a weak sharp minimum at infinity, i.e.,
		there exist $c>0$ and $R>0$ such that
		\begin{equation*}
			f(x)-f_*\geq c\,\mathrm{dist}\,(x,\Sol\eqref{problem-0}) \ \ \forall x\in\Omega\setminus\mathbb{B}_R.
		\end{equation*} 
		\item $f$ is coercive on $\Omega$.
	\end{enumerate}
\end{theorem}
\noindent{\it Proof\,}  Let $\tilde{f}\colon\R^n\to\overline{\R}$ be the function defined by $\tilde{f}(x):=(f+\delta_{\Omega})(x)$ for all $x\in\R^n$. Then $\mathrm{argmin}_{x\in\Omega} f(x)=\mathrm{argmin}_{x\in\R^n}\tilde{f}(x)$. By Proposition \ref{Pro-n-2}, one has
	\begin{equation*}
		\partial\tilde{f}(\infty)=\bigcup_{u\in \mathbb{S}}\partial\tilde{f}(\infty; u).
	\end{equation*}
	Furthermore, it follows from Proposition \ref{pro-3} that if $\partial\tilde{f}(\infty; u)$ is nonempty, then $u\in\Omega^\infty$. Indeed, if there exists $\xi\in \partial\tilde{f}(\infty; u)$, then we can find sequences $x_k\xrightarrow{u}\infty$ and $\xi_{k}\in\partial (f+\delta_{\Omega})(x_k)$ with $\xi_k\to\xi$. Hence, $x_k\in\dom f\cap\Omega$. This implies that $x_k\xrightarrow{\Omega, u}\infty$ and so $u\in\Omega^\infty$. Thus 
	\begin{equation*}
		\partial\tilde{f}(\infty)=\bigcup_{u\in \Omega^\infty\cap \mathbb{S}}\partial\tilde{f}(\infty; u).
	\end{equation*} 
	Since condition \eqref{equ-n-6} holds for all $u\in\Omega^\infty\cap \mathbb{S}$, one has
	\begin{equation*}
		\partial\tilde{f}(\infty; u)\subset \partial f(\infty; u)+N_\Omega(\infty; u) \ \ \forall u\in \Omega^\infty\cap \mathbb{S}.
	\end{equation*}	 
	Hence,
	\begin{equation*}
		\partial\tilde{f}(\infty)\subset\bigcup_{u\in \Omega^\infty\cap \mathbb{S}}[\partial f(\infty; u)+N_\Omega(\infty; u)].
	\end{equation*}
	This and condition \eqref{equ-n-5} imply that $0\notin \partial\tilde{f}(\infty)$ and conclusions of the theorem follow directly from \cite[Theorem 6.4]{Kim-Tung-Son-23}.  \qed  
 
\medskip
The following example demonstrates an application of Theorem \ref{Theorem3.3} and indicates that \cite[Theorem 6.4]{Kim-Tung-Son-23} is not applicable in this context.
\begin{example}\rm  Let $f\colon\mathbb{R}^2\to\mathbb{R}$,   $x\mapsto e^{-x_1}+(x_2-x_1)^2$ for all $x=(x_1, x_2)\in \mathbb{R}^2$ and let $\Omega:=\{x\in\mathbb{R}^2\,|\, x_1=0\}$. Then we have 
	$$\nabla f(x)=\big(-e^{-x_1} -2(x_2-x_1), 2(x_2-x_1)\big)\ \ \forall x\in\R^2.$$ 
	For each $r\in\R$ and $k\in \mathbb{N}$, we see that
	\begin{equation*}
		\nabla f(k-\tfrac{r}{2}, k)\to (r, -r)
	\end{equation*}	 
	as $k\to\infty$. Hence, $\partial f(\infty)=\{(r, -r)\,|\, r\in\R\}$ and so $0\in \partial f(\infty)+N_\Omega(\infty)$. This means that \cite[Theorem 6.4]{Kim-Tung-Son-23} cannot be employed for this example. We now use our  Theorem \ref{Theorem3.3} to show that $\mathrm{argmin}_{x\in\Omega}\,f(x)$ is nonempty and compact, and problem \eqref{problem-0} has a weak sharp minimum at infinity. Indeed, it is easy to see that  $\Omega^\infty\cap\mathbb{S}=\{(0, \pm 1)\}$. For $u=(0, 1)$, if $\xi=(\xi_1, \xi_2)\in \partial f(\infty; u)$, then there exists sequences $x_k=(x_{1k}, x_{2k})\xrightarrow{u} \infty$ with  $\nabla f(x_k)\to \xi$. This implies that $x_{2k}\to  +\infty$, $\tfrac{x_{1k}}{x_{2k}}\to 0$, and $x_{2k}-x_{1k}\to \tfrac{\xi_2}{2}$, which is impossible. Hence, $\partial f(\infty; u)=\emptyset$. Similarly, if $u=(0, -1)$, we have also $\partial f(\infty; u)=\emptyset$.  Thus 
	$$0\notin \partial f(\infty; u)+N_\Omega(\infty; u)\ \ \forall u\in \Omega^\infty\cap\mathbb{S},$$
	and the conclusions follow from Theorem \ref{Theorem3.3}.
\end{example}

\subsection{Error Bound at Infinity}
Let $g\colon\mathbb{R}^n\to\overline{\mathbb{R}}$ be an l.s.c. function and $\Omega$ be a nonempty and closed subset in $\mathbb{R}^n$. Consider the constraint set
\begin{equation}\label{constraint-set}
	S:=\{x\in \Omega\;|\; g(x)\leq 0\}.	
\end{equation}
Assume that $\Omega\cap\dom g$ is unbounded and $S$ is nonempty. 

\begin{definition}\rm 
	We say that the constraint set $S$ has an {\em error bound at infinity} if there exist $\alpha>0$ and $R>0$ such that
	\begin{equation}\label{equ-n-7}
		\mathrm{dist}\,(x; S)\leq \alpha [g(x)]_+ \ \ \forall x\in \Omega\setminus\mathbb{B}_R. 
	\end{equation}

	A following result presents a sufficient condition for the existence of error bounds at infinity for the constraint set $S$.
\end{definition}
\begin{theorem}\label{theorem-3.2} Let 
	$S$ be as in \eqref{constraint-set}. Assume that for all $u\in\Omega^\infty\cap\mathbb{S}$ the following conditions are satisfied: 
	\begin{equation*}\label{CQ-infinity}
		\partial^\infty g(\infty;u)\cap(-N_\Omega(\infty; u))=\{0\}
	\end{equation*}
	and 
	\begin{equation*}\label{regular-infinity}
		0\notin \partial g(\infty; u)+N_\Omega(\infty; u).
	\end{equation*}
	Then, $S$  has an error bound at infinity. 
\end{theorem}
\noindent{\it Proof\,}  Let $\tilde{g}\colon\R^n\to\overline{\R}$ be the function defined by $\tilde{g}(x):=g(x)+\delta_{\Omega}(x)$ for all $x\in\R^n$. Then 
	$$S=\{x\in\R^n\;|\; \tilde{g}(x)\leq 0\}$$ 
	and condition \eqref{equ-n-7} is equivalent to
	\begin{equation}\label{equ-n-8}
		\mathrm{dist}\,(x, S)\leq \alpha[\tilde{g}(x)]_+\ \ \forall x\in\R^n\setminus\mathbb{B}_R.
	\end{equation}
	Analysis similar to that in the proof of Theorem \ref{Theorem3.3} shows that $0\notin \partial \tilde{g}(\infty)$. This and \cite[Theorem 3.2]{Tuyen-24} imply that there exist $\alpha>0$ and $R>0$ satisfying \eqref{equ-n-8}. This means that $S$ has an error bound at infinity.  \qed
 
\medskip
The example below shows an application of Theorem \ref{theorem-3.2} and makes clear that \cite[Theorem 3.2]{Tuyen-24} cannot be employed in this context.
\begin{example}\rm Let $\Omega=\mathbb{R}_+$ and let $g\colon \mathbb{R} \to \mathbb{R}$ be the function defined by
	\begin{equation*}
		g(x)=
		\begin{cases}
			x, \ \ &\text{if}\ \ x\geq 0,
			\\
			e^x-1, \ \ &\text{if}\ \ x<0.
		\end{cases}
	\end{equation*}
	Then the constraint system \eqref{constraint-set} is
	\begin{equation*}
		S=\{x\in\Omega\;|\; g(x)\leq 0\}= \{0\}.
	\end{equation*}
	An easy computation shows that $\Omega^\infty\cap \mathbb{S}=\{1\}$,  and $\partial g(\infty; 1)=\{1\}$. Thus $0\notin \partial g(\infty; 1)+N_\Omega(\infty; 1)$ and so the constraint set $S$ has an error bound at infinity due to Theorem \ref{theorem-3.2}. On the other hand, we can see that $\partial g(\infty)=\{0, 1\}$ and  so $0\in \partial g(\infty)+N_\Omega(\infty)$. This implies that \cite[Theorem 3.2]{Tuyen-24} cannot be applied for this example. 
\end{example}

The following result is a corollary of Theorem  \ref{theorem-3.2} and Proposition \ref{pro411}. 

\begin{corollary}
	Let $S$ be a constraint set defined by
	\begin{equation*}\label{equa-21}
		S:=\{x\in \Omega\;:\; g_i(x)\leq 0, i\in I:=\{1, \ldots, m\}\}
	\end{equation*} 
	where $g_i\colon\R^n\to\R$, $i\in I$, are l.s.c. functions, $\Omega$ is an unbounded closed subset in $\R^n$ such that $\Omega\cap (\cap_{i\in I}\,\dom\, g_i)$ is unbounded. If for all $u\in\Omega^\infty\cap\mathbb{S}$, the following conditions hold
	\begin{equation*}\label{equa-11}
		\left[\xi_1+\ldots+\xi_m+v=0, \xi_i\in\partial^{\infty} g_i(\infty;u), \eta\in N_\Omega(\infty; u)\right] \Rightarrow \xi_i=\eta=0 \ \ \forall i\in I,
	\end{equation*} 
	and  
	\begin{equation*}\label{equa-12}
		\nexists \lambda\in\Delta_m \ \ \text{such that}\ \   0\in \sum_{i=1}^m\lambda_i\circ\partial g_i(\infty; u) +N_\Omega(\infty; u),
	\end{equation*}
	then, $S$ has an error bound at infinity, i.e., there exist $\alpha>0$ and $R>0$ such that
	\begin{equation*}\label{equa-14}
		d(x; S)\leq \alpha\sum_{i=1}^m[g_i(x)]_+ \ \ \forall x\in\Omega\setminus\mathbb{R}.
	\end{equation*}
\end{corollary}

\section{Conclusions}\label{Conclusion}
In this paper, we introduced and studied directional normal cones at infinity together with directional limiting and singular subdifferentials at infinity. We established calculus rules for these concepts and demonstrated their usefulness in nonsmooth optimization through several applications, including directional optimality conditions at infinity, coercivity of the objective function, compactness of the global solution set, and properties such as weak sharp minima and error bounds at infinity. Illustrative examples were provided to highlight the effectiveness of the proposed framework and to compare it with the existing results in \cite{Kim-Tung-Son-23,HA-Hung-25}.

We see the following natural directions for future developments in the theory of directionally generalized differentiation at infinity.

1. Motivated by the recent work of Kim et al. \cite{KPTT-23}, one should consider developing a version of directional coderivatives at infinity for set-valued mappings.

2. Using this aid together with the works \cite{Gfrerer-13,Gfrerer-132,Tung-Son}, one can establish criteria at infinity for the directional well-posedness properties of set-valued mappings.       


\section*{Funding} This research is funded by the Hanoi Pedagogical University 2 [grant number HPU2.2025-UT-04].

\section*{Disclosure statement} 
The authors declare that they have no conflict of interest.

\section*{Data availability} There is no data included in this paper.

\end{document}